\documentclass[final]{article}
\usepackage{amsfonts}
\usepackage{verbatim}

\usepackage{float}
\floatstyle{plain}
\newfloat{mfile}{tp}{mfi}[section]
\newfloat{myfig}{tp}{mfi}[section]
\floatname{mfile}{Listing}
\floatname{myfig}{Figure}

\def\figverbatiminput#1{{\footnotesize\verbatiminput{#1}}}

\usepackage{latexsym,amsmath,amssymb,graphics,amscd}
\usepackage{graphics,color}
\usepackage{graphicx}

\newtheorem{lemma}{Lemma}

\newcommand{\RR}{\mathbb{R}}
\newcommand{\hf}{{\textstyle{{\frac{1}{2}}}}}
\newcommand{\hhf}{{\scriptstyle{{\frac{1}{2}}}}}

\title{Deep Learning: An Introduction for Applied Mathematicians}

\author{
Catherine F. Higham\thanks{School of Computing Science, 
University of Glasgow, UK (\texttt{Catherine.Higham@glasgow.ac.uk}). This author was supported by the EPSRC UK Quantum Technology Programme under grant EP/M01326X/1.} 
\and
Desmond J. Higham\thanks{Department of Mathematics and Statistics, University of Strathclyde, UK (\texttt{d.j.higham@strath.ac.uk}).
This author was supported by the EPSRC/RCUK 
Digital Economy Programme
under grant  EP/M00158X/1.
}
}

\begin{document}

\maketitle

\begin{abstract}
Multilayered artificial neural networks are becoming a pervasive tool in a host of application fields. At the heart of this 
deep learning revolution are familiar concepts from applied and computational 
mathematics; notably, in calculus, approximation theory, optimization and linear algebra.
This article provides a very brief introduction to the basic ideas that underlie deep learning  from an applied mathematics perspective.
Our target audience includes postgraduate and final year undergraduate
students in mathematics who
are keen to learn about the area.
The article may also be useful for instructors in mathematics who wish to enliven their classes with references to 
the application of deep learning techniques.
We focus on three fundamental questions: 
what is  a deep neural network? 
how is a network trained? what is the stochastic gradient method?
We illustrate the ideas with a short MATLAB code that sets up and 
trains a network. We also show the use of state-of-the art software
on a large scale image classification problem.
We finish with references to the current literature. 
\end{abstract}




\section{Motivation}
\label{sec:mot}
 
Most of us have come across the phrase deep learning.
It relates to a set of tools that have become extremely popular 
in a vast range of application fields, from 
image recognition,
speech recognition
 and natural language processing 
 to targeted advertising and drug discovery.
The field has grown to the extent where sophisticated 
software packages are available in the public domain, many produced 
by high-profile technology companies. Chip manufacturers are also 
customizing their graphics processing units (GPUs)
for kernels at the heart of deep learning.

Whether or not its current level of attention is fully justified, 
deep learning is clearly a topic of interest to  
employers, and therefore to our students. 
Although there are many useful resources available, we feel that 
there is a niche for a brief treatment  
aimed at mathematical scientists. 
For a mathematics student, 
gaining some familiarity with deep learning can 
enhance employment prospects.
For mathematics educators, 
slipping \lq\lq Applications to 
Deep Learning\rq\rq\ into the syllabus 
of a class on 
   calculus,
   approximation theory,
   optimization,
   linear algebra,
   or scientific computing 
is a great way 
to attract students and maintain their interest in core topics.
The area is also ripe for independent project study.

There is no novel material in this article, and 
many topics are glossed over or omitted.
Our aim is to present some key ideas as 
clearly as possible while avoiding 
non-essential detail.
The treatment is aimed at readers with a background in 
mathematics who have completed a course in 
linear algebra and are familiar with partial differentiation. 
Some experience of scientific computing 
is also desirable.

To keep the material concrete, we 
list and walk through 
a short MATLAB code that 
illustrates the main algorithmic 
steps in setting up, training and applying an 
artificial neural network.
We also demonstrate 
the high-level use of state-of-the-art  
software on a larger scale problem.

Section~\ref{sec:annex}
introduces 
some key ideas by creating and 
training an artificial neural network
using a simple example.
Section~\ref{sec:gen}
sets up some 
useful notation and defines a 
general network.
Training a network, which involves the solution of an 
optimization problem, is the main computational challenge in this field.
In Section~\ref{sec:sg} we describe the 
stochastic gradient method, a variation 
of a traditional optimization technique that is designed to cope with very 
large scale sets of training data.
Section~\ref{sec:bp} 
explains how the 
partial derivatives 
needed for 
the 
stochastic gradient method 
can be computed efficiently using back propagation.
First-principles MATLAB code that illustrates 
these ideas is provided in 
section~\ref{sec:mat}.
A larger scale problem is treated in 
section~\ref{sec:image}. Here we 
make use of existing software.
Rather than repeatedly acknowledge 
work throughout the text, we have 
chosen to place the bulk of our citations in Section~\ref{sec:disc},
where pointers to the large and growing literature 
are given. In that section we also raise issues
that were not mentioned elsewhere, and highlight some 
current hot topics.

\section{Example of an Artificial Neural Network}
\label{sec:annex}
\begin{figure}
\begin{center}
\includegraphics[scale=0.4]{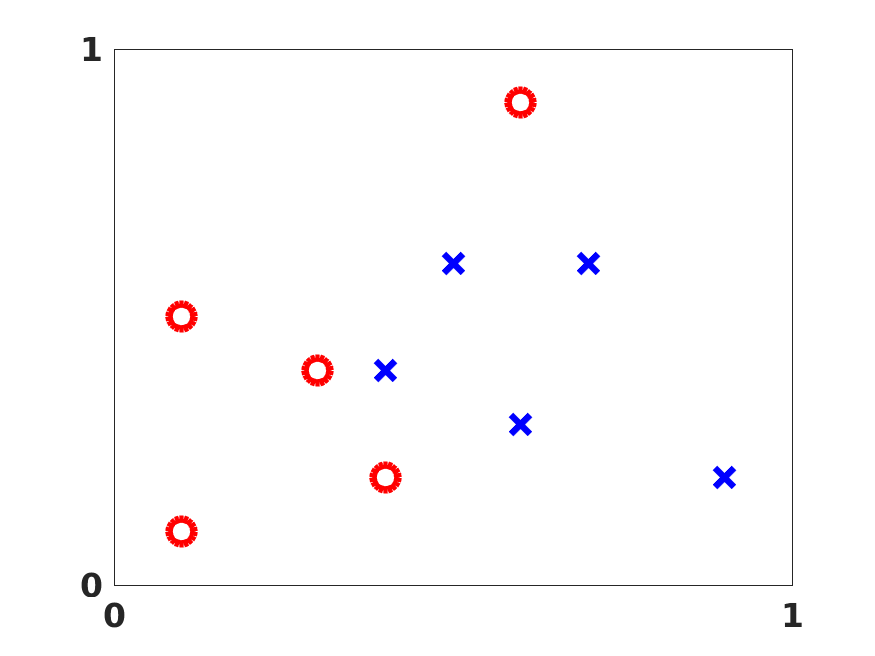}
\caption{
Labeled data points in $\RR^2$. 
Circles denote points in category A. 
Crosses denote points in category B. 
}
\label{Fig:pic_xy}
\end{center}
\end{figure}

This article takes a data fitting view of artificial neural networks.
To be concrete, consider the set of points shown in 
Figure~\ref{Fig:pic_xy}.
This shows \emph{labeled data}---some points are in category A, indicated by circles, and the
rest are in category B, indicated by crosses.
For example, the data may show oil drilling sites on a map, where category A  denotes 
a successful  outcome.
Can we use this data to categorize a newly proposed drilling site?
Our job is to construct a transformation that takes any point in $\RR^2$ 
and returns 
either a circle or a square.
Of course, there are many reasonable ways to construct such a transformation. 
The artificial neural network approach uses repeated application of a simple, nonlinear
function. 

\begin{figure}
\begin{center}
\includegraphics[scale=0.6]{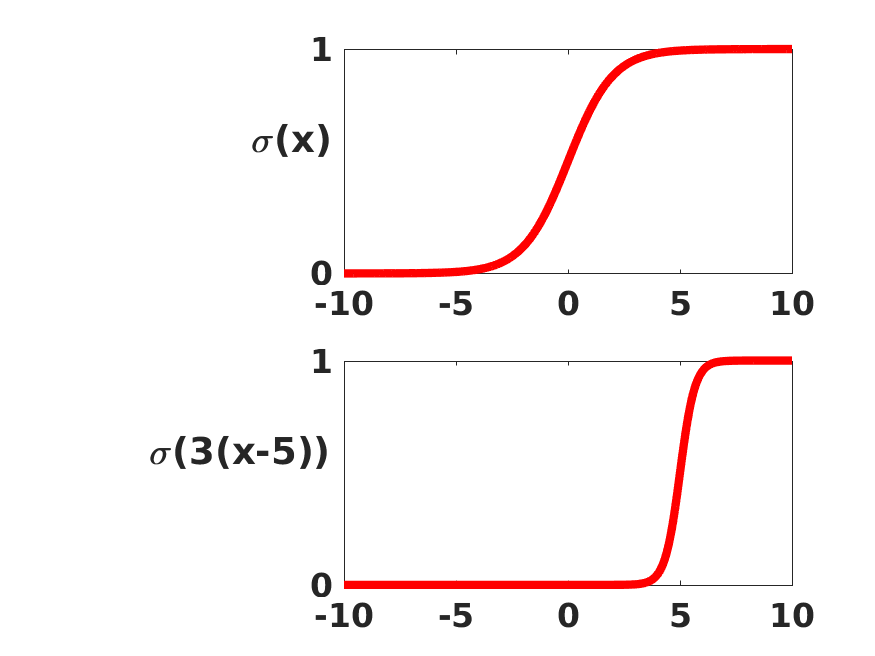}
\caption{
Upper: sigmoid function (\ref{eq:sig}). 
Lower: sigmoid with shifted and scaled input.
}
\label{Fig:sigpic}
\end{center}
\end{figure}

We will base our network on the sigmoid function 
\begin{equation}
 \sigma(x) = \frac{1}{ 1 + e^{-x}},
\label{eq:sig}
\end{equation}
which is illustrated in the upper half of Figure~\ref{Fig:sigpic} 
over the interval $-10 \le x \le 10$.
We may regard $\sigma(x)$ as a smoothed version of a step function, which itself mimics the  behavior of a neuron in the brain---firing
(giving output equal to one) if the input is large enough, and 
remaining inactive (giving output equal to zero)
otherwise. 
The sigmoid also has the convenient property that its 
derivative takes the simple form
\begin{equation}
 \sigma'(x) = \sigma(x) \left( 1 - \sigma(x) \right),
\label{eq:sigdiff}
\end{equation}
which is straightforward to verify.

The steepness and location of the transition in the sigmoid function may be altered by scaling and shifting the argument or, in the 
language of neural networks, by
\emph{weighting} and \emph{biasing} the input.
 The lower plot in 
Figure~\ref{Fig:sigpic} shows $\sigma\left(3(x-5)\right)$. 
The factor $3$ has sharpened the changeover 
and the shift $-5$ has altered its location.
To keep our notation managable, we need to interpret the sigmoid function in a vectorized sense.
For $z \in \RR^m$, 
$\sigma: \RR^m \to \RR^m$ is defined by applying the sigmoid function in the obvious 
componentwise manner, so that
\[
 \left( \sigma(z) \right)_i = \sigma( z_i ).
\] 
With this notation, we can set up layers of neurons.
In each layer, every neuron outputs a single real number, which is passed to every neuron in 
the next layer.
At the next layer, each neuron forms its own weighted combination of these values, adds its own 
bias, and applies the sigmoid function.
Introducing some mathematics, 
if the real numbers produced by the neurons in one 
layer are collected into a vector, $a$, then the vector of 
outputs from the next layer 
has the form 
\begin{equation}
      \sigma( W a + b).
\label{eq:swab}
\end{equation}
Here, $W$ is matrix and $b$ is a vector. 
We say that $W$ contains the \emph{weights} 
and $b$ contains the \emph{biases}.
The number 
of columns in $W$ matches the number of 
neurons that produced the vector $a$ at the previous layer. The number of 
rows 
in $W$ matches the number of 
neurons at the current layer.
The number of components in $b$ also 
matches the number of 
neurons at the current layer.
To emphasize the role of the $i$th neuron in 
(\ref{eq:swab}), we could pick out the $i$th component as 
\[
      \sigma \left( \sum_{j} w_{ij}  a_j + b_i \right),
\]
where the sum runs over all entries in $a$.
Throughout this article, we will be switching between the vectorized and componentwise viewpoints to strike a balance between 
clarity and brevity.

\begin{figure}
\begin{center}
\includegraphics[scale=0.6,angle=0]{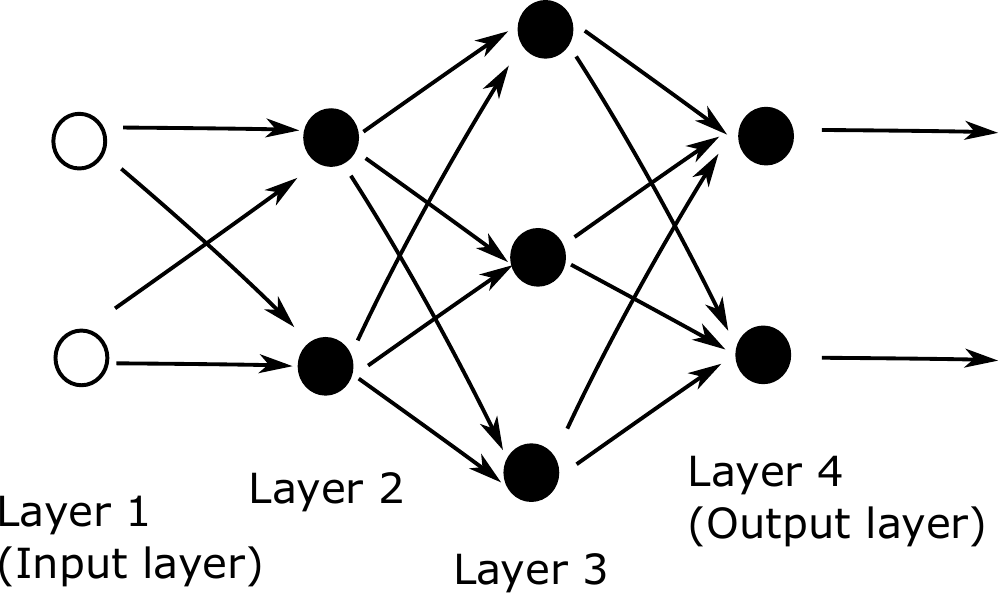}
\caption{
A network with four layers.}
\label{Fig:netpic4}
\end{center}
\end{figure}

In the next section, we introduce a full set of notation  
that allows us to define a general network. 
Before reaching that stage, we will give 
a specific example. 
Figure~\ref{Fig:netpic4} represents an artificial 
neural network with four layers.
We will apply this form of network
to the problem defined by Figure~\ref{Fig:pic_xy}.
For the network in 
Figure~\ref{Fig:netpic4} the first (input) layer is
represented by two circles. 
This is because our input data points have two components.
The second layer has two solid circles, 
indicating that two neurons are being employed.
The arrows from layer one to layer two indicate that both 
components of the input data 
are made available to the two neurons in layer two.
Since the input data has the form $x \in \RR^2$, the weights and biases for layer two may be 
represented by a matrix
$W^{[2]} \in \RR^{2 \times 2}$ and a vector
$b^{[2]} \in \RR^{2}$, respectively.
The output from layer two then has the form 
\[
  \sigma( W^{[2]} x + b^{[2]}) \in \RR^2.
\]
Layer three has three neurons, each receiving input in $\RR^2$. 
So the weights and biases
for 
layer three may be 
represented by a matrix
$W^{[3]} \in \RR^{3 \times 2}$ and a vector
$b^{[3]} \in \RR^{3}$, respectively.
The output from layer three then has the form 
\[
  \sigma\left( W^{[3]} 
     \sigma( W^{[2]} x + b^{[2]})
 + b^{[3]}\right) \in \RR^3.
\]
The fourth (output) layer has two neurons,  each receiving input in 
$\RR^3$. So the weights and biases
for this layer may be 
represented by a matrix
$W^{[4]} \in \RR^{2 \times 3}$ and a vector
$b^{[4]} \in \RR^{[2]}$, respectively.
The output from layer four, and hence from the overall network, 
has the form 
\begin{equation}
 F(x) = 
 \sigma\left( W^{[4]} 
  \sigma\left( W^{[3]} 
     \sigma( W^{[2]} x + b^{[2]})
 + b^{[3]}\right) 
  + b^{[4]}
 \right)
\in \RR^2.
\label{eq:all4}
\end{equation}

The expression (\ref{eq:all4}) defines a function $F: \RR^2 \to \RR^2$ in terms of its 
 $23$ parameters---the entries in the weight matrices and bias vectors.
Recall that our aim is to produce a classifier based on 
the data in 
Figure~\ref{Fig:pic_xy}.
We do this by optimizing over the parameters.
We will require $F(x)$ to be close to $[1,0]^T$ for data points 
in category A and 
  close to $[0,1]^T$ for data points in category B.
Then, given a new point $x \in \RR^2$, it would be reasonable to classify it according to the largest 
component  of $F(x)$; that is, category A if $F_1(x) > F_2(x)$ and 
category B if  $F_1(x) < F_2(x)$,
with some rule to break ties.
This requirement on $F$ may be specified through a \emph{cost function}.
Denoting the ten data points in 
Figure~\ref{Fig:pic_xy}
by 
$\{ x^{\{ i \}} \}_{i=1}^{10}$, we use 
$y( x^{\{ i \}} )$ for the target output; that is, 
\begin{equation}
y( x^{\{ i \}} ) = \left\{ \begin{array}{l} 
                                    \left[ \begin{array}{c}
                                                 1 \\
                                                 0
                                            \end{array}
                                      \right]
                                          \quad \mathrm{if~}
x^{\{ i \}} \mathrm{~is~in~category~A}, \\
                  ~~ \\ 
                                    \left[ \begin{array}{c}
                                                 0 \\
                                                1 
                                            \end{array}
                                      \right]
                                          \quad \mathrm{if~}
x^{\{ i \}} \mathrm{~is~in~category~B}. 
                      \end{array}
                \right.
\label{eq:ytarget}
\end{equation}
Our cost function then takes the form 
\begin{equation}
\mathrm{Cost}
   \left(W^{[2]},W^{[3]},W^{[4]}, b^{[2]}, b^{[3]}, b^{[4]}\right)
     = 
      \frac{1}{10}
             \sum_{i=1}^{10}
                \hf \| 
              y( x^{\{ i \}} ) - F( x^{\{ i \}} )
                               \|_2^2.
 \label{eq:C10}
\end{equation}
Here, the factor $\hhf$ is included for convenience; 
it simplifies matters when we start differentiating.
We emphasize that 
$\mathrm{Cost}$ is a function of the weights and biases---the data points 
are fixed.
The form in (\ref{eq:C10}), where discrepancy is 
measured by averaging the squared Euclidean norm 
over the data points, is
often refered to as a \emph{quadratic cost function}.
In the language of optimization, $\mathrm{Cost}$  is our
\emph{objective function}. 

Choosing the weights and biases in a way that minimizes 
the cost function is refered to  
as
\emph{training} the network.
We note that, 
in principle, rescaling an objective function 
does not change an optimization problem. 
We should arrive at the same minimizer 
if we change 
$\mathrm{Cost}$ 
to, for example, 
$100 \, \mathrm{Cost}$ or
$\mathrm{Cost}/30$.
So the factors 
$1/10$ and $1/2$ in
(\ref{eq:C10})
should have no effect on the outcome.

\begin{figure}
\begin{center}
\includegraphics[scale=0.4]{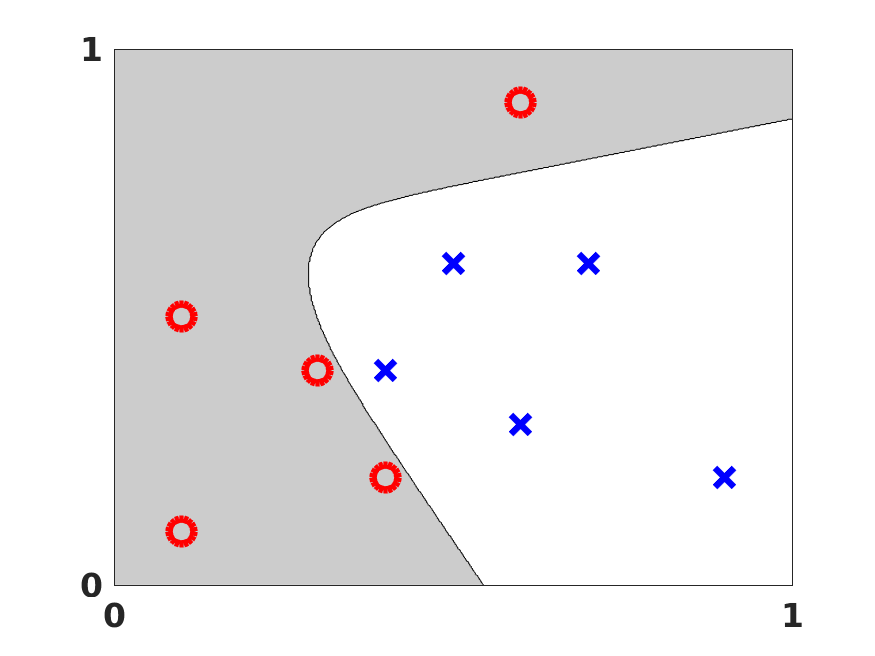}
\caption{
Visualization of output from an artificial neural network 
applied to the data in Figure~\ref{Fig:pic_xy}.
}
\label{Fig:pic_bdy}
\end{center}
\end{figure}

For the data in 
Figure~\ref{Fig:pic_xy}, we used 
the MATLAB optimization toolbox to minimize the 
cost function 
(\ref{eq:C10}) over the $23$ parameters defining  
   $W^{[2]}$, $W^{[3]}$, $W^{[4]}$, $b^{[2]}$, $b^{[3]}$ and $b^{[4]}$.
More precisely, we used 
the nonlinear least-squares solver 
\texttt{lsqnonlin}.
For the trained network,  
Figure~\ref{Fig:pic_bdy}
 shows the boundary where $F_1(x) > F_2(x)$. So, with this approach, any point
in the shaded region would be assigned to category A and 
any point
in the unshaded region to category B.

\begin{figure}
\begin{center}
\includegraphics[scale=0.4]{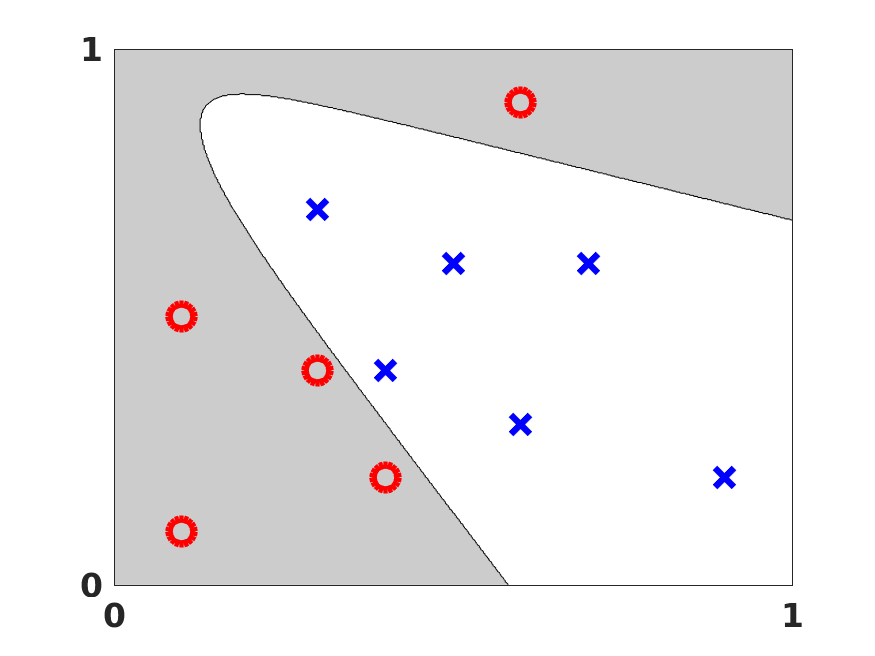}
\caption{
Repeat of the experiment in  
Figure~\ref{Fig:pic_bdy} with an additional data point.
}
\label{Fig:pic_bdy2}
\end{center}
\end{figure}

Figure~\ref{Fig:pic_bdy2} shows how the network responds to 
additional training data.
Here we added one further category B point,  
indicated by the extra cross at $(0.3,0.7)$, 
and re-ran the optimization routine.

The example illustrated in 
Figure~\ref{Fig:pic_bdy}
is small-scale by the standards of today's 
deep learning tools. However, the underlying 
optimization problem, minimizing a non-convex 
objective function over $23$ variables, is fundamentally difficult.
We cannot exhaustively search over a 23 dimensional parameter 
space, and we cannot guarantee to find the global minimum of a 
non-convex function.
Indeed, some experimentation with the location of the data points 
in 
Figure~\ref{Fig:pic_bdy}
and with the choice of initial guess for the weights and biases 
makes it clear that 
\texttt{lsqnonlin}, with its default settings, 
cannot always find an acceptable solution.
This motivates the 
material in sections~\ref{sec:sg}
and \ref{sec:bp}, where we focus on the 
optimization problem.

\section{The General Set-up}
\label{sec:gen}
The four layer example in section~\ref{sec:annex} 
introduced the idea of neurons, represented by the sigmoid function, 
acting in layers. 
At a general layer,
each neuron receives the same input---one real value from every neuron at the previous layer---and produces one real value, which is passed to every neuron at the next layer.
There are two exceptional layers. At the input layer, there is no ``previous layer'' and each 
neuron receives the input vector. At the output layer, there is no ``next layer'' and these neurons provide the overall output.
The layers in between these two are called \emph{hidden layers}. There is no 
 special meaning to this phrase; it simply indicates that these neurons are performing intermediate calculations.
Deep learning is a loosely-defined term which implies that many hidden layers are being used.
 
We now  
spell out the general form of the notation from section~\ref{sec:annex}.
We suppose that the network has $L$ layers, with layers
$1$ and $L$ being the input and output layers, respectively.
Suppose that layer $l$, for $l = 1, 2, 3, \ldots, L$ contains 
$n_l$ neurons.
So $n_1$ is the dimension of the input data.
Overall, the network maps from 
$\RR^{n_1}$ to $\RR^{n_L}$.
We use 
$W^{[l]} \in \RR^{n_l \times n_{l-1}}$ to denote 
the matrix of weights at layer $l$.
More precisely,  $w^{[l]}_{jk}$
is the weight that neuron $j$ at layer $l$ applies to the output 
from neuron $k$ at layer $l-1$.
Similarly, $b^{[l]} \in \RR^{n_l}$ is the vector 
of biases for layer $l$, so neuron 
$j$ at layer $l$ uses the bias $b^{[l]}_j$.

\begin{figure}
\begin{center}
\includegraphics[scale=0.6]{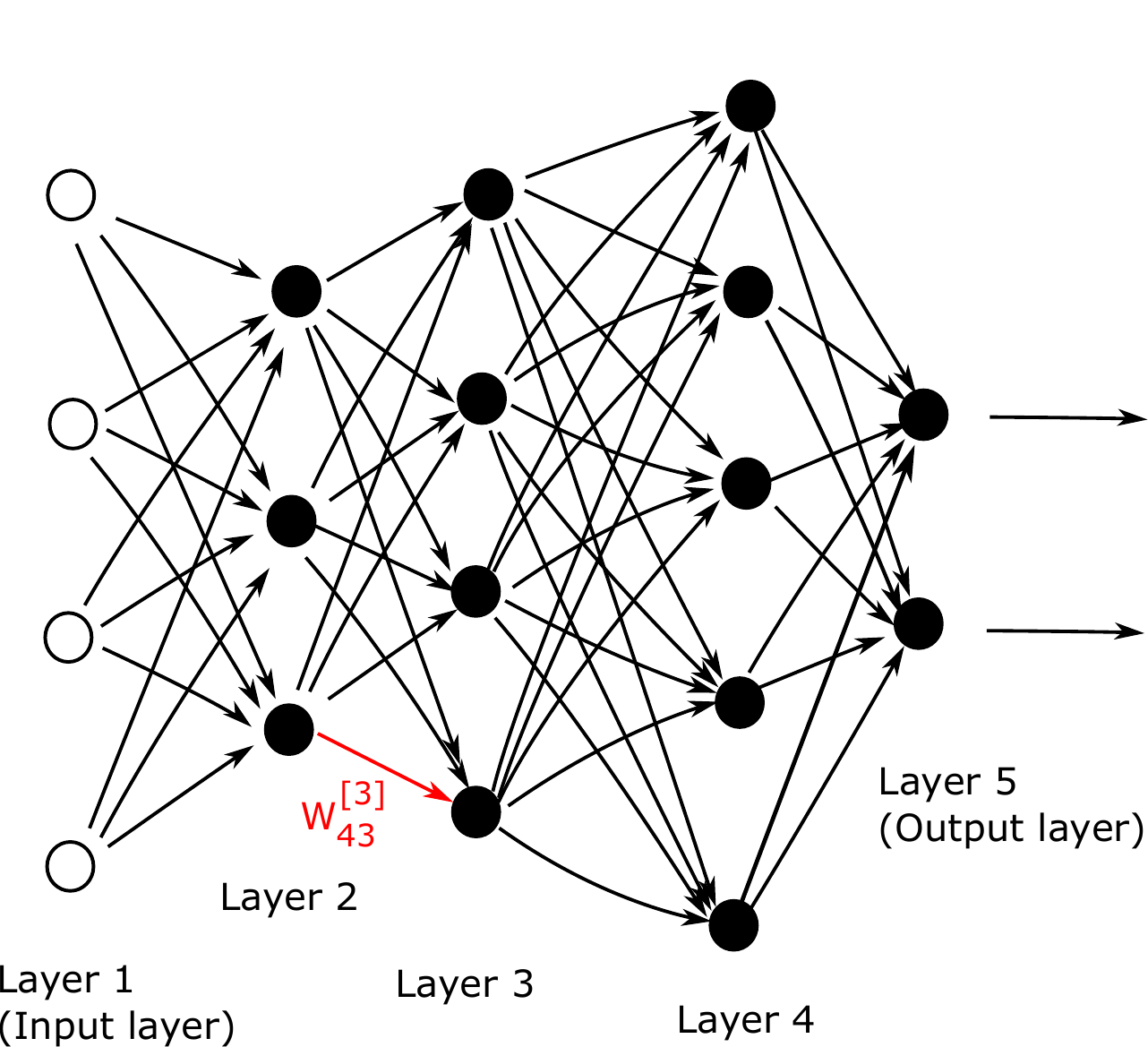}
\caption{
A network with five layers.
The edge corresponding to the weight 
$w_{43}^{[3]}$ 
is highlighted.
The output from neuron number 3 at layer 2 is weighted 
by the factor 
$w_{43}^{[3]}$ 
when it is fed into neuron number 4 at layer 3.
}
\label{Fig:netpic5}
\end{center}
\end{figure}

In Fig~\ref{Fig:netpic5}
we give an example with $L = 5$ layers. 
Here, 
$n_1 = 4$, 
$n_2 = 3$, 
$n_3 = 4$, 
$n_4 = 5$ and 
$n_5 = 2$, so 
 $W^{[2]} \in \RR^{3 \times 4}$,
 $W^{[3]} \in \RR^{4 \times 3}$,
 $W^{[4]} \in \RR^{5 \times 4}$,
 $W^{[5]} \in \RR^{2 \times 5}$,
$b^{[2]} \in \RR^{3}$,
$b^{[3]} \in \RR^{4}$,
$b^{[4]} \in \RR^{5}$
and $b^{[5]} \in \RR^{2}$.

Given an input $x \in \RR^{n_1}$,  
we may then neatly summarize the action of the network  
by letting 
$a_j^{[l]}$ denote the output, or \emph{activation},  
from neuron $j$ at layer $l$.
So, we have 
\begin{eqnarray}
 a^{[1]} &=& x \in \RR^{n_1},
\label{eq:a1def}
\\
 a^{[l]} &=& \sigma \left( W^{[l]} a^{[l-1]} + b^{[l]} \right)
 \in \RR^{n_l}, \quad \mathrm{for~} l = 2, 3, \ldots, L.
\label{eq:aldef}
\end{eqnarray} 
It should be clear that 
(\ref{eq:a1def}) and (\ref{eq:aldef}) amount to an algorithm for feeding
the input forward through the network in order to produce an output 
$a^{[L]} \in \RR^{n_L}$. At the end of 
section~\ref{sec:bp} this 
algorithm appears within a pseudocode description of an 
approach for training a network.

Now suppose we have $N$ pieces of data, or 
\emph{training points}, in $\RR^{n_1}$, 
$\{ x^{\{i\}} \}_{i=1}^{N}$, 
for which there are given target outputs
$\{ y( x^{\{i\}} ) \}_{i=1}^{N}$ in $ \RR^{n_L}$.
Generalizing 
(\ref{eq:C10}), the quadratic cost function
that we wish to minimize has the form 
\begin{equation}
\mathrm{Cost} 
     = 
      \frac{1}{N}
             \sum_{i=1}^{N}
                \hf \| 
              y( x^{\{ i \}} ) - a^{[L]}( x^{\{ i \}} )
                               \|_2^2,
 \label{eq:costfun}
\end{equation}
where, to keep notation under control, we have not explicitly  
indicated that $\mathrm{Cost}$ is a function of all the weights and biases.

\section{Stochastic Gradient}
\label{sec:sg}
We saw in the previous two sections that training a network corresponds to 
choosing the parameters, that is, the weights and biases, that minimize the cost function.
The weights and biases take the form of matrices and vectors, but 
at this stage 
it is convenient 
to imagine them stored as a single vector that we call $p$. 
The example in Figure~\ref{Fig:netpic4} has a total of 23 weights and biases. So, in that case, $p \in \RR^{23}$.
Generally, we will suppose  $p \in \RR^{s}$, and write the 
cost function 
in 
(\ref{eq:costfun}) as $\mathrm{Cost} (p)$ to emphasize its dependence on
the parameters.
 So 
 $\mathrm{Cost} : \RR^{s} \to \RR$.

We now introduce a classical method in optimization that is often referred to as
\emph{steepest descent}
or
\emph{gradient descent}.
The method proceeds iteratively, computing a sequence of vectors in  $\RR^{s}$
with the aim of converging to a vector that minimizes 
the cost function.
Suppose that our current vector is $p$.
How should we choose a perturbation, $\Delta p$, so that the next vector, 
$p + \Delta p$, represents an improvement? 
If $\Delta p$ is small, then ignoring terms of order $\| \, \Delta p \, \|^2$, a Taylor series
expansion gives
\begin{equation}
\mathrm{Cost} ( p + \Delta p) \approx  \mathrm{Cost} (p) + \sum_{r = 1}^{s}  
          \frac{ \partial \, \mathrm{Cost}(p) } {\partial \, p_r}    \Delta p_r.
\label{eq:Costexp}
\end{equation}
Here $\partial \, \mathrm{Cost}(p)  / \partial \, p_r $
denotes the partial derivative of the cost function with respect to the $r$th parameter.
For convenience, we will let $ \nabla \mathrm{Cost}(p)  \in \RR^{s}$ denote the vector of 
partial derivatives, known as the \emph{gradient}, so that 
\[
  \left( \nabla \mathrm{Cost}(p) \right)_{r} =  \frac{ \partial \, \mathrm{Cost}(p) } {\partial \, p_r}.
\]
Then (\ref{eq:Costexp}) becomes 
\begin{equation}
\mathrm{Cost} ( p + \Delta p) \approx  \mathrm{Cost} (p) 
        + 
                \nabla \mathrm{Cost}(p)^T \Delta p.
\label{eq:Costexpb}
\end{equation}
Our aim is to reduce the value of the cost function. 
The relation (\ref{eq:Costexpb}) motivates the idea of choosing  
$\Delta p$ to make $\nabla \mathrm{Cost}(p)^T \Delta p $ as negative as possible.
We can address this problem via the Cauchy--Schwarz inequality, which states that for 
any $f,g \in    \RR^{s}$, we have 
$| f^T g | \le \| \, f \, \|_2 \| \, g \, \|_2$.
So the most negative that $f^T g$ can be is 
$ -\| \, f \, \|_2 \| \, g \, \|_2$, which happens when 
$f = -g$.
Hence, based on (\ref{eq:Costexpb}),  we should choose 
$\Delta p$ to lie in the direction $-\nabla \mathrm{Cost}(p)$.
Keeping in mind that (\ref{eq:Costexpb}) is an approximation that is relevant only for small
$\Delta p$, we will limit ourselves to a small step in that direction.
This leads to the update 
\begin{equation}
  p \to p - \eta \nabla \mathrm{Cost}(p).
\label{eq:gd}
\end{equation}
Here $\eta$ is small stepsize that, in this context, is known as the 
\emph{learning rate}.
This equation defines the 
steepest descent method. We choose an initial vector and 
iterate with 
(\ref{eq:gd}) until some stopping criterion has been met, or 
until the number of iterations has exceeded the computational budget.

Our cost function 
 (\ref{eq:costfun}) involves  
a sum 
of 
individual terms that runs 
over the training data. 
It follows that the partial derivative 
  $\nabla \mathrm{Cost}(p)$ 
is a sum 
over the training data 
of individual  
partial derivatives.  
More precisely, let 
\begin{equation}
C_{x^{ \{ i \} }} = 
                \hf \| 
              y( x^{\{ i \}} ) - a^{[L]}( x^{\{ i \}} )
                               \|_2^2.
\label{eq:Cxi}
\end{equation}
Then, from 
(\ref{eq:costfun}), 
\begin{equation}
\nabla 
\mathrm{Cost} (p) 
     = 
      \frac{1}{N}
             \sum_{i=1}^{N}
                   \nabla 
C_{x^{ \{ i \}}} (p).
 \label{eq:pdcostfun}
\end{equation}

When we have a large number of parameters and a large number 
of training points, computing 
the gradient vector 
(\ref{eq:pdcostfun}) at every iteration of the steepest descent 
method 
(\ref{eq:gd}) can be prohibitively expensive.
A much cheaper alternative is to replace the mean of 
the individual gradients over all training points by the gradient at a
single, randomly chosen, training point. 
This leads to the simplest form of what is called 
the \emph{stochastic gradient} method. 
A single step may be summarized as 
\begin{enumerate}
 \item Choose an integer $i$ uniformly at random from 
$\{1,2,3,\ldots,N\}$.
 \item Update 
\begin{equation}
 p \to p 
-  
 \eta 
                   \nabla C_{ x^{ \{ i \}} } (p).
\label{eq:sgdef}
\end{equation}
\end{enumerate}
In words, at each step, 
the stochastic gradient method 
uses one randomly chosen training point to 
represent the full training set.
As the iteration proceeds, the 
method 
sees more training points. So there is some hope that 
this dramatic reduction in cost-per-iteration will be worthwhile
overall. 
We note that, even for very small $\eta$,  
the update
(\ref{eq:sgdef}) is not guaranteed to reduce the 
overall cost function---we have traded the mean for a single 
sample. 
Hence, although the phrase \emph{stochastic gradient descent}
is widely used, we prefer to 
use stochastic gradient. 

The version of the stochastic gradient method that we introduced in 
(\ref{eq:sgdef}) 
is the simplest from a large range of possibilities. 
In particular, 
the index $i$ in 
(\ref{eq:sgdef})
was chosen by 
sampling \emph{with replacement}---after using a training point, 
it is returned to the training set and is just as likely
as any other point to be 
chosen at the next step.
An alternative is to sample without 
replacement; that is, to cycle through each of the 
$N$ training points in a random order.
Performing $N$ steps in this manner, refered to as completing an  
\emph{epoch}, 
may be summarized as follows:
\begin{enumerate}
 \item Shuffle the integers 
       $\{1,2,3,\ldots,N\}$ into a new order, 
       $\{k_1,k_2,k_3,\ldots,k_N\}$.
 \item for $i = 1$ upto $N$, update 
\begin{equation}
p \to p 
- 
 \eta 
                   \nabla C_{ x^{ \{ k_i \}} } (p).
\label{eq:sgdefrep}
\end{equation}
\end{enumerate}

If we regard the 
stochastic gradient method 
as approximating the mean over all 
training points in 
 (\ref{eq:pdcostfun}) by a single sample, 
then it is natural to 
consider a compromise where we use 
a small sample average. For some $m \ll N$ 
we could take steps of the following form. 
\begin{enumerate}
 \item Choose $m$ integers, $k_1, k_2, \ldots, k_m$,  
uniformly at random from  
$\{1,2,3,\ldots,N\}$.
 \item Update 
\begin{equation}
 p \to p 
- 
 \eta 
   \frac{1}{m} 
     \sum_{i = 1}^{m}
                   \nabla C_{ x^{ \{ k_i \}} } (p).
\label{eq:sgdefbatch}
\end{equation}
\end{enumerate}
In this iteration, the set 
$
          \{ 
            x^{ \{ k_i \}}
            \}_{i=1}^{n}
$ 
is known as a \emph{mini-batch}.
There is a \emph{without replacement} alternative where, 
assuming $N = Km$ for some $K$, we split the training set 
randomly into $K$ distinct mini-batches and cycle through them.
  
Because the stochastic gradient method is usually 
implemented within the context of a very large scale 
computation, algorithmic choices such as 
mini-batch size and the form of randomization 
are often driven by the 
requirements of high performance computing 
architectures.
Also, it is, of course, possible to vary these choices, along with 
others, such as the learning rate, 
dynamically as the training progresses in an attempt to 
accelerate convergence. 

Section~\ref{sec:mat}
describes a simple MATLAB code that uses a vanilla 
stochastic gradient method.
In section~\ref{sec:image} we use a 
state-of-the-art implementation and 
section~\ref{sec:disc}
has pointers to the current literature.

\section{Back Propagation}
\label{sec:bp}
We are now in a position to apply the 
stochastic gradient method in order to train an 
artificial neural network.
So we switch from the general vector of parameters, $p$,  
used in section~\ref{sec:sg} to the 
entries in the weight matrices and bias vectors.
Our task is to 
compute partial derivatives of the cost function with respect to 
each $w^{[l]}_{jk}$ and 
$b^{[l]}_j$.
We saw that 
the idea behind the stochastic 
gradient method is to exploit the 
structure of the 
cost function: because 
(\ref{eq:costfun}) 
is a linear combination of individual terms that runs 
over the training data the same is true
of its partial derivatives. 
We therefore focus our attention on computing those 
individual partial derivatives. 

Hence, for a fixed training 
point we regard 
$
   C_{x^{\{ i \}}} 
$ 
in (\ref{eq:Cxi})
as a function of the weights and biases.
So we may drop the dependence on 
$x^{\{ i \}}$
and simply write 
\begin{equation}
   C = \hf \| y - a^{[L]} \|_2^2.
\label{eq:Cx}
\end{equation}
We recall from 
(\ref{eq:aldef}) 
that $a^{[L]}$ is the output from the 
artificial neural network. The dependence of 
$C$ on the weights and biases arises only through $a^{[L]}$.

To derive worthwhile expressions for the partial derivatives, 
it is useful to introduce two further sets of variables.
First we 
let 
\begin{equation}
 z^{[l]} = W^{[l]} a^{[l-1]} + b^{[l]} \in \RR^{n_l}, \qquad
 \mathrm{for~} l = 2,3,\ldots,L.
\label{eq:azl}
\end{equation}
We refer to 
$z^{[l]}_{j}$ as the \emph{weighted input} for neuron $j$ at layer $l$.
The fundamental relation 
(\ref{eq:aldef}) 
that propagates information through the network  
may then be written 
\begin{equation}
   a^{[l]} = \sigma  \left( z^{[l]} \right), \quad \mathrm{~for~} l = 2,3,\ldots,L.
\label{eq:aldefb}
\end{equation}
Second, we let 
$
 \delta^{[l]}
 \in \RR^{n_l}$ 
 be defined by 
\begin{equation}
 \delta^{[l]}_j = \frac{ \partial \, C}{\partial \, z^{[l]}_j}, 
  \qquad 
    \mathrm{~for~} 1 \le j \le n_l \mathrm{~~~and~~~}  2 \le l \le L.
\label{eq:dellj}
\end{equation}
This expression, which is often called the \emph{error}
in the $j$th neuron at layer $l$,
is an
intermediate quantity 
that is useful 
both for 
analysis and computation.
However, we point out that this useage of the 
term error is somewhat ambiguous.
At a general, hidden layer, it is not clear how much to 
\lq\lq blame\rq\rq\ each neuron for discrepancies in 
the final output.
Also, at the output layer, $L$, the expression 
(\ref{eq:dellj}) does not quantify those discrepancies directly.
The idea of referring to 
$
 \delta^{[l]}_j 
$
in (\ref{eq:dellj}) as an error 
seems to have arisen because 
the cost function can only be at a minimum if all partial 
derivatives are zero, so 
$
 \delta^{[l]}_j  = 0 
$
is a useful goal.
As we mention later in this section, it may be more helpful to 
keep in mind that 
$
 \delta^{[l]}_j   
$
measures the sensitivity of 
the cost function to the weighted input for neuron $j$ at layer $l$.

At this stage we also need to define the 
Hadamard, or
componentwise,  
product of two vectors. If $x,y \in \RR^n$, then
$x \circ y \in \RR^n$ is defined by $(x \circ y)_{i} = 
x_i y_i$. In words, the Hadamard product is formed by
pairwise multiplication of the corresponding components.

With this notation, the following results are a consequence of the 
chain rule.

\begin{lemma}
\label{lem:bp}
We have 
\begin{eqnarray}
  \delta^{[L]} &=&  \sigma'(z^{[L]} )  \circ (a^L - y),
\label{eq:delL}
\\
  \delta^{[l]} &=&  \sigma'(z^{[l]} ) 
                       \circ 
   ( W^{[l+1]} )^T \delta^{[l+1]},
 \qquad 
 \qquad 
 \qquad 
\mathrm{~for~} 2 \le l \le L-1,
\label{eq:dell}
\\
 \frac{ \partial \, C} {\partial \, b^{[l]}_j} &=& \delta^{[l]}_j,
 \qquad 
 \qquad 
 \qquad 
 \qquad 
 \qquad 
 \qquad 
 \qquad 
 \quad 
\mathrm{~for~} 2 \le l \le L,
\label{eq:pCb}
\\
 \frac{ \partial \, C} {\partial \, w^{[l]}_{jk}} &=& 
\delta^{[l]}_{j}
      a^{[l-1]}_{k}.
 \qquad 
 \qquad 
 \qquad 
 \qquad 
 \qquad 
 \qquad 
 \quad 
  \mathrm{~for~} 2 \le l \le L.
\label{eq:pCw}
\end{eqnarray}
\end{lemma}

\textbf{Proof}
We begin by proving 
(\ref{eq:delL}).
The relation 
(\ref{eq:aldefb}) 
with $l = L$ shows that 
$z^{[L]}_j$ and $a^{[L]}_j$ are connected by 
$a^{[L]} = \sigma \left( z^{[L]} \right)$, and hence 
\[
  \frac{ \partial \, a^{[L]}_j }{ \partial \, z^{[L]}_j } 
= \sigma'(z^{[L]}_j).
\]
Also, from (\ref{eq:Cx}), 
\[
   \frac{ \partial \, C}{\partial \, a^{[L]}_j } =
   \frac{ \partial }{\partial \, a^{[L]}_j }
      \hf \sum_{k = 1}^{n_L} (y_k - a^{[L]}_k)^2 
   =   -( y_j - a^{[L]}_j).
\]
So, using the chain rule,
\[
  \delta^{[L]}_j = \frac{ \partial \, C}{\partial \, z^{[L]}_j} =
  \frac{ \partial \, C}{\partial \, a^{[L]}_j} 
  \frac{ \partial \, a^{[L]}_j}{\partial \, z^{[L]}_j} =
   (a^{[L]}_j - y_j) 
    \sigma'(z^{[L]}_j),
\]
which 
is the componentwise form of 
(\ref{eq:delL}).

To show 
(\ref{eq:dell}),
we use the chain rule to convert from $z_j^{[l]}$ to 
$\{ z^{[l+1]}_{k} \}_{k=1}^{n_{l+1}}$.
Applying the chain rule, and using the definition
(\ref{eq:dellj}),
\begin{equation}
   \delta^{[l]}_j = \frac{\partial \, C} {\partial z^{[l]}_j } 
                        = \sum_{k = 1}^{n_{l+1}} 
                                        \frac{\partial \, C} {\partial z^{[l+1]}_k } 
                                           \frac{\partial \, z^{[l+1]}_k} {\partial  z^{[l]}_j}
                  =  \sum_{k = 1}^{n_{l+1}} 
                                      \delta^{[l+1]}_k 
                                           \frac{\partial \, z^{[l+1]}_k} {\partial  z^{[l]}_j}. 
\label{eq:pfc}
\end{equation} 
Now, from (\ref{eq:azl}) 
we know that 
$z^{[l+1]}_k$ and $z^{[l]}_j$ are connected via
\[
  z^{[l+1]}_k = \sum_{s=1}^{n_l} w^{[l+1]}_{ks} \sigma\left( z^{[l]}_s \right) + b^{[l+1]}_k.
\]
Hence, 
\[
 \frac{\partial \, z^{[l+1]}_k} {\partial  z^{[l]}_j} 
   = w^{[l+1]}_{kj} \sigma'\left( z^{[l]}_j \right).
\] 
In (\ref{eq:pfc}) this gives 
\[
  \delta^{[l]}_j = 
       \sum_{k = 1}^{n_{l+1}} 
                                      \delta^{[l+1]}_k w^{[l+1]}_{kj} \sigma'\left( z^{[l]}_j \right),
\]
which may be rearranged as
\[
     \delta^{[l]}_j =   \sigma'\left( z^{[l]}_j \right) \left( (W^{[l+1]})^T  \delta^{[l+1]} \right)_j.
\]
This is the componentwise form of 
(\ref{eq:dell}). 

To show
(\ref{eq:pCb}), 
 we note from 
(\ref{eq:azl}) 
and 
(\ref{eq:aldefb})
that 
$z^{[l]}_j$ is connected to $b^{[l]}_j$ by 
\[
   z^{[l]}_j = \left( W^{[l]} \sigma\left( z^{[l-1]} \right) \right)_j + b^{[l]}_j.
\]
Since $z^{[l-1]}$ does not depend on $b^{[l]}_j$,
we find that 
\[
  \frac{ \partial \, z_j^{[l]} } { \partial \, b^{[l]}_j} = 1.
\]
Then, from the chain rule,
\[
 \frac{ \partial \, C }{ \partial \, b^{[l]}_j }    
     =
    \frac{ \partial \, C }{ \partial \, z^{[l]}_j }   
    \frac{ \partial  \, z^{[l]}_j}{ \partial \, b^{[l]}_j }   
   = 
      \frac{ \partial \, C }{ \partial \, z^{[l]}_j }
  = \delta^{[l]}_j,
\]
using the 
definition 
(\ref{eq:dellj}). This gives 
(\ref{eq:pCb}). 

Finally, to obtain 
(\ref{eq:pCw}) 
 we start with the componentwise version of 
(\ref{eq:azl}), 
\[
  z^{[l]}_j = \sum_{k = 1}^{n_{l-1}} w^{[l]}_{jk} a^{[l-1]}_k + b^{[l]}_j,
\]
which gives
\begin{equation}
   \frac{ \partial \, z^{[l]}_j } { \partial w^{[l]}_{jk} }
    = a^{[l-1]}_{k},   \qquad \mathrm{independently~of~} j,
\label{eq:pfd}
\end{equation}
and
\begin{equation}
   \frac{ \partial \, z^{[l]}_s } { \partial w^{[l]}_{jk} }
    = 0,   \qquad \mathrm{for~} s \neq j.
\label{eq:pfe}
\end{equation}
In words, 
(\ref{eq:pfd}) 
and 
(\ref{eq:pfe}) 
follow because the $j$th neuron at layer $l$ uses the weights from only the $j$th row of 
$W^{[l]}$, and applies these weights linearly.
Then, from the chain rule, 
(\ref{eq:pfd}) and 
(\ref{eq:pfe}) give 
\[
   \frac{ \partial \, C} { \partial \, w^{[l]}_{jk} }  = 
  \sum_{s=1}^{n_l}
     \frac{ \partial \, C} { \partial z^{[l]}_s} 
    \frac{ \partial \, z^{[l]}_s }{ \partial \, w^{[l]}_{jk} } 
     = 
     \frac{ \partial \, C} { \partial z^{[l]}_j } 
     \frac{ \partial \, z^{[l]}_j }{ \partial \, w^{[l]}_{jk} }
     = 
     \frac{ \partial \, C} { \partial z^{[l]}_j } 
     a^{[l-1]}_k 
   =
     \delta^{[l]}_j    a^{[l-1]}_k ,
\] 
where the last step used the definition of $\delta^{[l]}_j $ in 
(\ref{eq:dellj}).
This completes the proof.
\quad $\blacksquare$

There are many aspects of Lemma~\ref{lem:bp}
that deserve our attention.
We recall from (\ref{eq:a1def}), 
(\ref{eq:azl})
and
(\ref{eq:aldefb}) 
that the output $a^{[L]}$ can be evaluated from a 
\emph{forward 
pass} 
through the network, computing
$a^{[1]}$, 
$z^{[2]}$, 
$a^{[2]}$, 
$z^{[3]}$, 
\ldots, 
$a^{[L]}$  
in order.
Having done this, 
we see from
(\ref{eq:delL})
that $\delta^{[L]}$ is immediately available.
Then, from 
(\ref{eq:dell}),
$\delta^{[L-1]}$,
$\delta^{[L-2]}$,
\ldots, 
$\delta^{[2]}$ may be computed 
in a 
\emph{backward 
pass}.
From 
(\ref{eq:pCb}) and 
(\ref{eq:pCw}), we then have access to the partial derivatives.
Computing gradients in this way is known 
as \emph{back propagation}.

To gain further understanding of the back propagation formulas
(\ref{eq:pCb}) 
and 
(\ref{eq:pCw}) 
in 
Lemma~\ref{lem:bp}, it is useful to recall the fundamental 
definition of a partial derivative.
The quantity 
 $
 \partial \, C / \partial \, w^{[l]}_{jk}
$ 
measures
how $C$ changes when we make a small perturbation to 
$
w^{[l]}_{jk}
$.
For illustration, 
Figure~\ref{Fig:netpic5} highlights the weight 
$
w^{[3]}_{43}
$.
It is clear that a change in this weight has no effect on the 
output of previous layers. So 
to work out 
$
 \partial \, C / \partial \, w^{[3]}_{43} 
$
we do not need to know about 
partial derivatives at previous layers.
It should, however, be possible to express   
$
 \partial \, C / \partial \, w^{[3]}_{43} 
$
in terms of 
partial derivatives at subsequent layers.
More precisely, 
the activation feeding into the $4$th neuron on layer $3$ 
is $z^{[3]}_{4}$, and, by definition, 
$\delta^{[3]}_{4}$ measures the sensitivity of 
$C$ with respect to this input.
Feeding in to this neuron we have 
$ 
w^{[3]}_{43}
a^{[2]}_3 
+ \mathrm{constant}$, so 
it makes 
sense that 
\[
\frac{ \partial \, C} { \partial \, w^{[3]}_{43} } 
   = \delta^{[3]}_{4} a^{[2]}_3.
\]
Similarly, in terms of the bias, 
$ b^{[3]}_4 
+ \mathrm{constant}$ is feeding in to the neuron, which explains why 
\[
\frac{ \partial \, C} { \partial \, b^{[3]}_{4} } 
 = 
   \delta^{[3]}_{4} \times 1.
\]

We may avoid the Hadamard product notation 
in
(\ref{eq:delL})
and 
(\ref{eq:dell})
by 
introducing diagonal matrices.
Let $D^{[l]} \in \RR^{n_l \times n_l}$ denote 
the diagonal matrix 
with $(i,i)$ entry given by 
$\sigma'(z^{[l]}_i)$.
Then we see that
  $\delta^{[L]} = D^{[L]} ( a^{[L]} - y)$ and 
  $\delta^{[l]} = D^{[l]} (W^{[l+1]})^T \delta^{[l+1]}$.
We could expand this out as 
\[
     \delta^{[l]} = 
  D^{[l]} (W^{[l+1]})^T
  D^{[l+1]} (W^{[l+2]})^T \cdots 
  D^{[L-1]} (W^{[L]})^T 
  D^{[L]}
   ( a^{[L]} - y).
\]
We also recall from 
(\ref{eq:sigdiff}) that $\sigma'(z)$ is trivial to compute. 

The relation 
(\ref{eq:pCb})
shows that 
$\delta^{[l]}$ corresponds precisely to the 
gradient of the cost function with respect to the biases at layer  
$l$.
If we regard 
$\partial C / \partial w^{[l]}_{jk}$ as defining the $(j,k)$
component in a matrix of partial derivatives at layer $l$, 
then 
(\ref{eq:pCw})
shows this matrix to be the \emph{outer product} 
$
\delta^{[l]} {a^{[l-1]}}^T \in \RR^{n_{l} \times n_{l-1}}
$.


Putting this together, we may write the following pseudocode for an 
algorithm that trains a network using a fixed number, 
$\mathrm{Niter}$, of stochastic gradient iterations. 
For simplicity, we consider the basic version 
(\ref{eq:sgdef})
where single samples are chosen with replacement.
For each 
training point,  
we perform a forward pass through the network in order to evaluate the 
activations, 
weighted inputs
and 
overall output 
$a^{[L]}$.
Then 
we 
perform a backward pass to compute the errors and updates.

\smallskip 

  \begin{description}
  \item For $\mathrm{counter} = 1$ upto $\mathrm{Niter}$
      \begin{description}
       \item Choose an integer $k$ uniformly at random from 
         $\{1,2,3,\ldots,N\}$
        \item $x^{ \{ k \}}$ is current training data point
        \item $a^{[1]}  = x^{ \{ k \}}$
        \item For $l = 2$ upto $L$
         \begin{description}
        \item $z^{[l]} = W^{[l]} a^{[l-1]} + b^{[l]}$ 
        \item $a^{[l]} = \sigma \left( z^{[l]} \right)$
        \item $D^{[l]} = \mathrm{diag} \left( \sigma'(  z^{[l]}  ) \right)$ 
        \end{description}
        \item end 
        \end{description}
        \begin{description}
      \item $\delta^{[L]} = D^{[L]} \left(a^{[L]} - y( x^{ \{ k \}} ) \right)$
      \item For $l = L-1$ downto $2$
        \begin{description}
      \item $\delta^{[l]} = D^{[l]} (W^{[l+1]})^T \delta^{[l+1]}$
     \end{description}
          \item  end
     \end{description}
     \begin{description}
        \item For $l = L$ downto $2$
      \begin{description}
        \item  $ W^{[l]} \to  W^{[l]} - \eta  \,
       \delta^{[l]} 
         {a^{[l-1]}}^T$  
            \item 
        $b^{[l]} \to  b^{[l]} - \eta   \,
       \delta^{[l]}$
     \end{description}
   \item end
   \end{description}
   \item end
   \end{description}

\section{Full MATLAB Example} \label{sec:mat}

We now give a concrete illustration involving back propagation and the
stochastic
gradient method.
Listing~\ref{List.netbp}
shows how a network of the form shown in Figure~\ref{Fig:netpic4}
may be used on the data in Figure~\ref{Fig:pic_xy}.
We note that this MATLAB code has been written for clarity  
and brevity, rather than efficiency or elegance.
In particular, we have \lq\lq hardwired\rq\rq\ the number of layers and
iterated 
through the forward and backward passes line by line.
(Because the weights and biases do not have the the same dimension in
each layer, it is 
not convenient to store them in a three-dimensional array.
We could use a cell array or structure array, 
\cite{MG3},
and
then implement the forward and backward passes in \texttt{for} loops.
However, this approach produced a less readable code, 
and violated our self-imposed one page limit.)

The function \texttt{netbp} in Listing~\ref{List.netbp}
contains the nested function \texttt{cost}, which evaluates 
a scaled version of 
$\mathrm{Cost}$ 
in (\ref{eq:C10}). 
Because this function is nested, it has access to the variables in the
main function, notably the training data.
We point out that the nested function \texttt{cost} is not used directly in the forward and
backward passes. It 
is called at each iteration of the stochastic gradient method so that we
can monitor the 
progress of the training.  

Listing~\ref{List.activate} shows the function \texttt{activate}, 
used by 
\texttt{netbp}, 
which applies the sigmoid function in vectorized form.

At the start of \texttt{netbp} we set up the training data and target 
$y$
values, as defined in 
(\ref{eq:ytarget}).
We then initialize all weights and biases using the normal pseudorandom
number generator
\texttt{randn}. For simplicity, we set a constant learning rate
\texttt{eta = 0.05} and 
perform a fixed number of iterations \texttt{Niter = 1e6}.

We use the the basic stochastic gradient iteration summarized at the end
of Section~\ref{sec:bp}. 
Here, the command \texttt{randi(10)} returns a uniformly and
independently chosen integer between $1$ and $10$.

Having stored the value of the cost function at each iteration, we use
the 
\texttt{semilogy}
command 
to visualize the progress of the iteration.

\begin{figure}
\begin{center}
\includegraphics[scale=0.4]{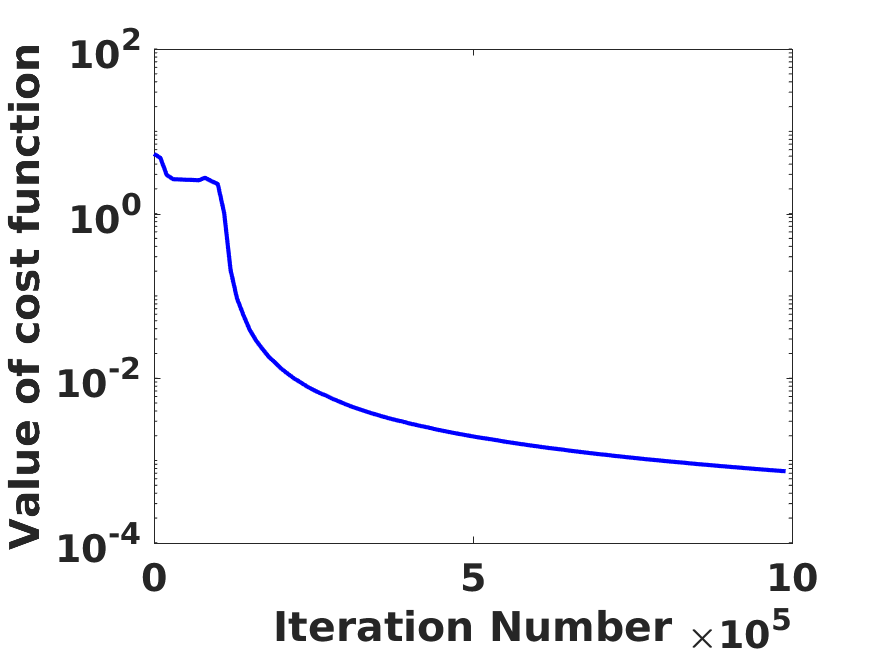}
\caption{
Vertical axis shows a scaled value of the 
cost function  
(\ref{eq:C10}).
Horizontal axis shows the iteration number.
Here we used the stochastic gradient method to train a network of the 
form shown in Figure~\ref{Fig:netpic4}
on the data in Figure~\ref{Fig:pic_xy}.
The resulting classification function is illustrated in 
Figure~\ref{Fig:pic_bdy_bp}.
}
\label{Fig:pic_bp_cost}
\end{center}
\end{figure}

In this experiment, our initial guess for the weights and biases 
produced
a cost function value of 
$5.3$.
After ${10}^6$ stochastic gradient steps this was 
reduced to $7.3 \times {10}^{-4}$.
Figure~\ref{Fig:pic_bp_cost}
shows the \texttt{semilogy} plot, and we see that the decay is not
consistent---the cost undergoes a flat period towards the start of the
process.
After this plateau,
we found that the cost decayed at a very slow linear rate---the 
ratio between successive values was typically within around
${10}^{-6}$ of unity.

\begin{figure}
\begin{center}
\includegraphics[scale=0.4]{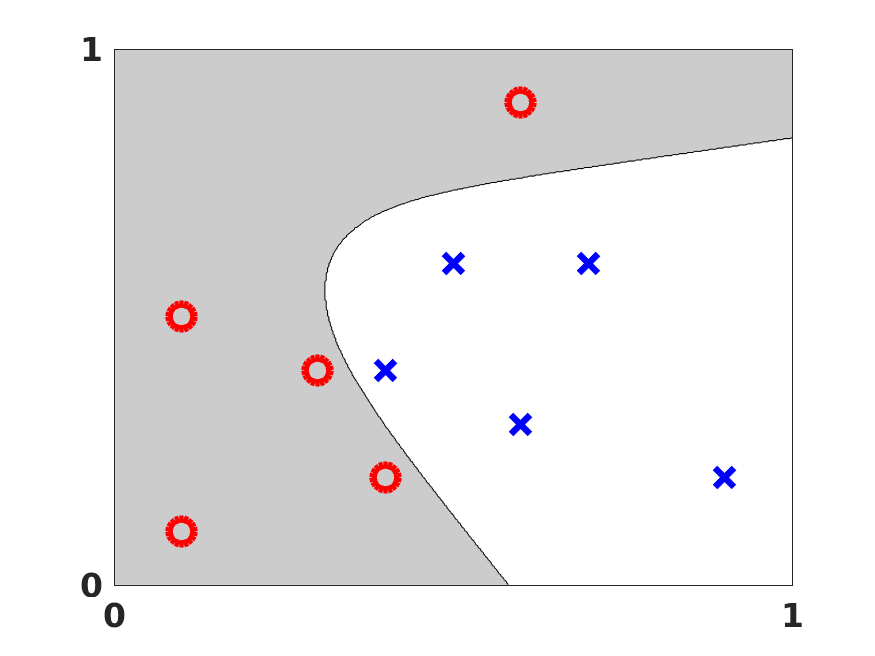}
\caption{
Visualization of output from an artificial neural network 
applied to the data in Figure~\ref{Fig:pic_xy}.
Here we trained the network 
using the stochastic gradient method with 
back propagation---behaviour of cost function is shown in 
Figure~\ref{Fig:pic_bp_cost}. The same optimization problem was
solved with the
\texttt{lsqnonlin} 
routine from MATLAB in order to produce 
Figure~\ref{Fig:pic_bdy}.
}
\label{Fig:pic_bdy_bp}
\end{center}
\end{figure}

\begin{figure}
\begin{center}
\includegraphics[scale=0.4]{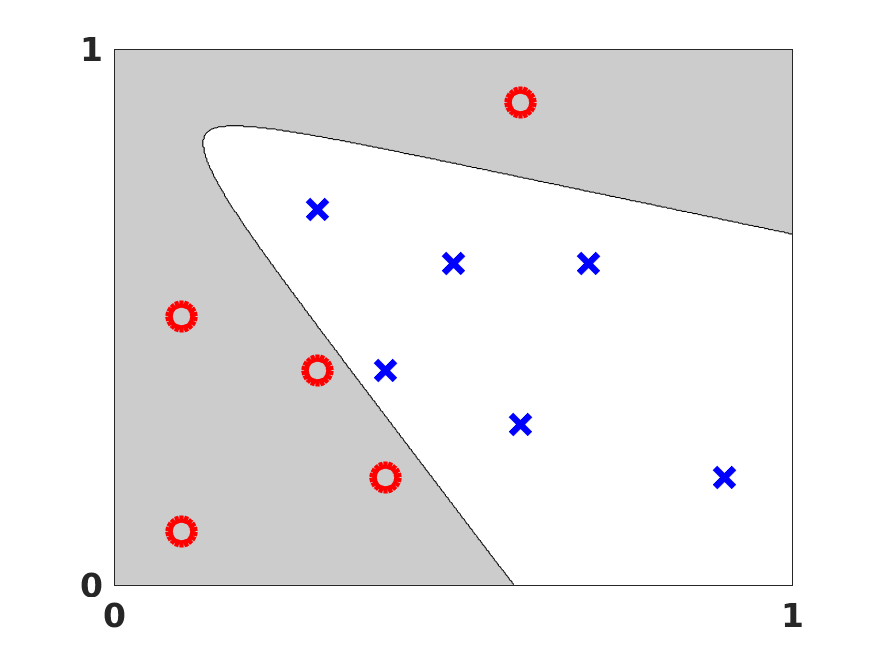}
\caption{
Visualization of output from an artificial neural network 
applied to the data in Figure~\ref{Fig:pic_xy} with an additional 
data point.
Here we trained the network 
using the stochastic gradient method with 
back propagation. The same optimization problem was
solved with the
\texttt{lsqnonlin} 
routine from MATLAB in order to produce 
Figure~\ref{Fig:pic_bdy2}.
}
\label{Fig:pic_bdy_bp2}
\end{center}
\end{figure}

An extended version of \texttt{netbp} can be found 
in the supplementary material.
This version has the extra graphics commands that make 
Figure~\ref{Fig:pic_bp_cost}
more
readable. It also 
takes the trained network
and produces Figure~\ref{Fig:pic_bdy_bp}.
This plot 
shows how the trained network carves up the input space.
Eagle-eyed readers will spot that the solution 
in 
Figure~\ref{Fig:pic_bdy_bp}.
differs 
slightly from the version
in 
Figure~\ref{Fig:pic_bdy}, 
where the same optimization problem was tackled by the
nonlinear least-squares 
solver \texttt{lsqnonlin}.
In 
Figure~\ref{Fig:pic_bdy_bp2}
we show the corresponding result when an extra data point
is added; this can be compared with 
Figure~\ref{Fig:pic_bdy2}.

\begin{mfile}

\caption{M-file \texttt{netbp.m}.}
\label{List.netbp}

\begin{center}
\figverbatiminput{netbp.m}
\end{center}

\end{mfile}

\begin{mfile}

\caption{M-file \texttt{activate.m}.}
\label{List.activate}

\begin{center}
\figverbatiminput{activate.m}
\end{center}

\end{mfile}

\section{Image Classification Example} \label{sec:image}

We now move on to a more realistic task, 
which allows us to demonstrate the power of the deep learning 
approach. 
We make use of the 
{\sc{MatConvNet}}
toolbox
\cite{Vedaldi15}, which is designed to 
offer key deep learning building blocks as 
simple MATLAB commands.
So 
{\sc{MatConvNet}} 
is an excellent environment for prototyping 
and for educational use.
Support for GPUs 
 also makes {\sc{MatConvNet}} efficient for large scale 
computations, and pre-trained networks may be downloaded for 
immediate use.

Applying 
{\sc{MatConvNet}} 
on a large scale problem 
also gives us the opportunity to outline  
further concepts that are relevant to practical computation. 
These are introduced in the next few subsections, before we  
apply them to the image classification exercise.

\subsection{Convolutional Neural Networks} 
\label{subsec:cnn}
{\sc{MatConvNet}}
uses a special class of
artificial neural networks known as 
Convolutional Neural Networks (CNNs), which have become a 
standard tool in computer vision applications.
To motivate CNNs, we note that 
the general framework described in 
section~\ref{sec:gen} does not scale well in the case of 
digital image data.
Consider a color image made up of 
200 by 200 pixels, each with a red, green and blue component.
This corresponds to an input vector 
of dimension $n_1 = 200 \times 200 \times 3 = 120,000$, and hence
a weight matrix $W^{[2]}$ at level 2 that has 
$120,000$ columns. 
If we allow general weights and biases, then this approach
is clearly infeasible. 
CNNs get around this issue by constraining the values 
that are allowed in the weight matrices and bias vectors.
Rather than a single full-sized linear transformation,  CNNs 
repeatedly apply a small-scale linear kernel, or filter, across 
portions of their input data.
In effect,  the weight matrices 
used by CNNs are extremely sparse and highly structured.

To understand why this approach might be useful, 
consider 
premultiplying an input vector in $\RR^6$ by the matrix 
\begin{equation}
\left[
 \begin{array}{rrrrrr}
               1 & -1 &     &  &  & \\
                 &  1 & -1  &  &  & \\
                 &    &  1  & -1 & & \\ 
                 &    &    & 1 & -1 & \\ 
                 &    &    &        &   1    & -1 
 \end{array}
\right]
\in \RR^{5 \times 6}.
\label{eq:toe}
\end{equation}
This 
produces a vector in $\RR^{5}$ made up of differences between 
neighboring values. 
In this case we are using a filter $[1, -1]$ and a \emph{stride} of 
length one---the filter advances by one place after each use. 
Appropriate generalizations of 
this matrix to the case of input vectors arising from 2D images 
can be used to detect \emph{edges} in an 
image---returning a large absolute 
value when there is an 
abrupt change in neighboring pixel values.
Moving a filter across an image can also reveal other  
features, for example, 
particular types of curves or 
blotches of the same color.  
So, having specified a filter size and stride length, we can allow the training process
to learn the weights in the filter as a means to extract useful structure. 

The word \lq\lq convolutional\rq\rq\ arises because the linear 
transformations involved may be written in the form of a 
convolution. In the 1D case, the convolution of the vector 
$x \in \RR^{p}$ with the filter  
$g_{1-p},  g_{2-p}, \ldots, g_{p-2}, g_{p-1}$ has $k$th component
given by 
\[
 y_{k} = \sum_{n = 1}^{p}   x_n g_{k-n}.
\]
The example in (\ref{eq:toe}) corresponds to a filter with 
$g_{0} = 1$, $g_{-1} =  -1$ and all other $g_k = 0$.
In the case
\begin{equation}
\left[
 \begin{array}{r}
    y_1 \\ 
    y_2 \\
    y_3 \\
    y_4   
 \end{array}
\right]
=
\left[
 \begin{array}{rrrrrrrrrr}
               a & b & c   & d  & & & & & &   \\
               & &  a & b & c   & d  & & & &   \\
               & & & &  a & b & c   & d  & &    \\
                & & & & & &  a & b & c   & d  
 \end{array}
\right]
\left[
 \begin{array}{c}
    x_1 \\ 
    x_2 \\
    x_3 \\
    x_4 \\
    x_5  \\
     x_6 \\ 
    x_7 \\
    x_8 \\
    x_9  \\
    0
 \end{array}
\right],
\label{eq:toep}
\end{equation}
we are applying a filter with four weights, $a$, $b$, $c$, and $d$, using a stride length of two.
Because the dimension of the input vector $x$ is not compatible with the filter length, we have 
\emph{padded} with an extra zero value.

In practice, 
image data is typically regarded as a three dimensional tensor: each pixel has two spatial coordinates and a red/green/blue value.
With this viewpoint, the filter takes the form of a  small tensor that is 
successsively applied to patches of the input tensor and 
the corresponding convolution operation is multi-dimensional. 
From a computational perspective, a key   
benefit of CNNs is that the matrix-vector products
involved in the forward and backward passes through the network 
can be computed extremely efficiently using fast 
transform techniques.

A convolutional layer is often followed by a 
\emph{pooling} layer, which reduces dimension by mapping 
small regions of pixels into single numbers.
For example, 
when 
these small regions are taken to be 
squares of four 
neigboring pixels in a 2D image,  
a 
\emph{max pooling} 
or 
\emph{average pooling} 
layer 
replaces each set of four
by their maximum or average value, respectively.

\subsection{Avoiding Overfitting}
\label{subsec:overf}
\emph{Overfitting} occurs when a trained network 
performs very accurately on the given data, but 
cannot generalize well to new data.
Loosely, this means that the fitting process has 
focussed too heavily on the unimportant and
unrepresentative \lq\lq noise\rq\rq\ in the 
given data.  
Many ways to combat overfitting have been suggested, some of 
which can be used together. 

One useful technique is to split the given data into two distinct 
groups.
\begin{itemize}
 \item \emph{Training data} is used in the 
       definition of the cost function that 
       defines the optimization problem. 
       Hence this data drives the process that iteratively 
       updates the weights.
 \item \emph{Validation data} is not used in the optimization
       process---it has no effect on the way that the weights are 
       updated from step to step.
       We use the validation data only to judge the performance of 
       the current network. At each step of the optimization, 
       we can evaluate the cost function 
       corresponding to the  
       validation data. This measures how well the 
       current set of trained weights performs on unseen data.
\end{itemize}
Intuitively, overfitting corresponds to the situation 
where the optimization process is driving down its 
cost function (giving a better fit to 
the training data), but the cost function 
for the validation error is no longer decreasing 
(so the performance on unseen data is not improving).
It is therefore reasonable to terminate 
the training at a stage where no improvement is seen on the 
validation data.

Another popular approach to tackle overfitting is to randomly and
independently remove neurons during the
training phase. For example, at each step of the stochastic gradient method, we could
delete each neuron with probability $p$ and train on the
remaining network.
At the end of the process, because the weights and biases 
were produced on these smaller networks, we could multiply each 
by a factor of $p$ for use on the full-sized network.
This technique, known as \emph{dropout},  
has the intuitive interpretation that we are
constructing an average over many trained networks,
with such a consensus
being more reliable than any individual.

\subsection{Activation and Cost Functions}
\label{subsec:ac}
In Sections~\ref{sec:annex} to \ref{sec:mat} we used 
activation functions of sigmoid form
(\ref{eq:sig}) and a quadratic cost function 
(\ref{eq:costfun}). There are many other widely used choices, 
and their relative performance is application-specific.
In our image classification setting it is common to use a 
\emph{rectified linear unit}, or ReLU, 
\begin{equation}
 \sigma(x) = \left\{ \begin{array}{c}
               0, \quad \mathrm{for~} x \le 0, \\
               x, \quad \mathrm{for~} x > 0, 
               \end{array}
                 \right.
\label{eq:relu}
\end{equation}
as the activation. 

In the case where our training data
$ \{ x^{ \{ i \} } \}_{i=1}^{N}$ comes from $K$ labeled categories,
let $l_i \in \{1,2,\ldots,K\}$ be the given label for data point $  x^{ \{ i \} } $.
As an alternative to the 
 quadratic cost function   
(\ref{eq:costfun}),
we could use a \emph{softmax log loss} approach, as follows.
Let the output  $a^{[L]}(x^{ \{ i \} } ) =: v^{ \{i\}}$
from the network take the form of a vector
in $ \RR^K$ such that the $j$th component is large when the image is believed to 
be from category $j$.
The \emph{softmax} operation
\[
  (v^{ \{i\}})_{s} \mapsto \frac{ e^{ v^{ \{i\}}_{s}  }} { \sum_{j = 1}^{K} e^{ v^{ \{i\}}_{j} } }.
\]
boosts the large components and produces a vector of positive weights 
summing to unity, which may 
be interpreted as probabilties.
	Our aim is now to force the softmax value for training point $ x^{ \{ i \} } $ to be as close to 
unity as possible in component $l_i$, 
which corresponds to the correct label.
Using a logarithmic rather than quadratic measure of error, we arrive at the cost function  
\begin{equation}
   - \sum_{i = 1}^{N}  
   \log \left(
   \frac{ e^{ v^{ \{i\}}_{l_i}  }} { \sum_{j = 1}^{K} e^{ v^{ \{i\}}_{j} }}
 \right).
\label{eq:smlcost}
\end{equation}

\subsection{Image Classification Experiment}
\label{subsec:it}

We now show results for a supervised learning task in
image classification.
To do this, we rely on 
the codes \verb5cnn_cifar.m5
and 
\verb5cnn_cifar_init.m5
that are available via 
the 
{\sc{MatConvNet}} website. 
We made only minor edits, including 
some changes that allowed us 
to test the use of dropout.
Hence, 
in particular, we are using 
the network architecture and parameter choices 
from those codes.
We refer to the 
{\sc{MatConvNet}} documentation and tutorial material for 
the fine details, and focus here on some of the 
bigger picture issues.

We consider a set of images, each of which 
falls into exactly one of the following ten categories: 
airplane, automobile, bird, cat, deer, dog, frog, horse, ship, truck.
We use labeled data from the 
freely available 
CIFAR-10 collection 
\cite{Krizhevsky09learningmultiple}.
The images are small, having $32$ by $32$ pixels, each with a red, green, and blue component.
So one piece of training data consists of $32 \times 32 \times 3 = 3,072$ values.   
We use a training set of 50,000  images, and 
use 10,000 more images as our validation set.
Having completed the optimization and trained the network, 
we then judge its performance on a fresh collection of 
10,000 test images, with 1,000 from each category.

\begin{figure}
\begin{center}
\includegraphics[scale=0.35]{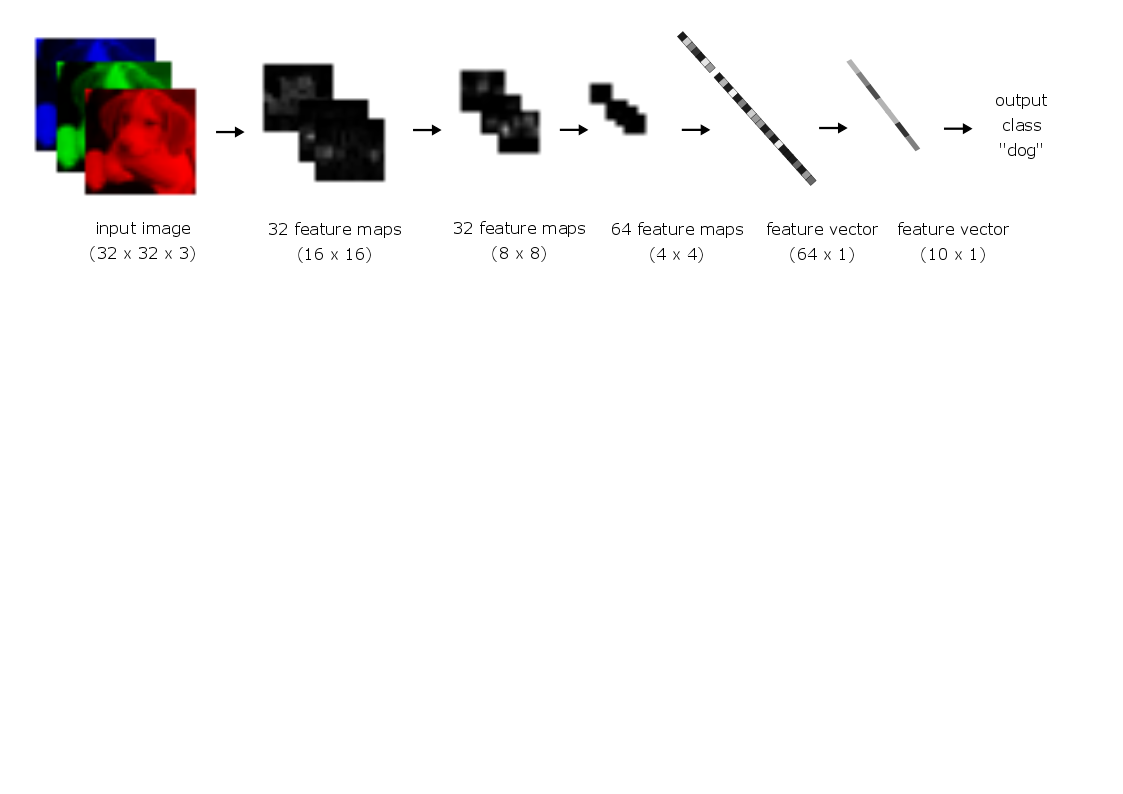}
\vspace{-2.5in}
\caption{
Overview of the CNN used for the image classification task. 
}
\label{Fig:netblocks}
\end{center}
\end{figure}

Following the 
architecture used in the 
relevant 
{\sc{MatConvNet}}  
codes, 
we set up a network whose layers are divided into five blocks as follows.
Here we 
describe the dimemsions of the inputs/outputs and
weights
in  
compact tensor notation. 
(Of course, the tensors could be stretched into sparse vectors and matrices in order 
to fit in with the general framework of sections~\ref{sec:annex}
to 
\ref{sec:mat}. But we feel that the tensor notation is natural in this context, and it is
consistent with the  {\sc{MatConvNet}} syntax.)
 \begin{description}
\item[Block 1] consists of a 
a convolution layer 
followed by a pooling layer and activation.
This
converts the original 
$32 \times 32 \times 3 $ input into dimension 
 $ 16 \times 16 \times 32$.
In more detail, the convolutional layer uses $5 \times 5$ filters that also scan across the $3$ color channels. There are $32$ different filters, so overall the weights
can be represented in a $5 \times 5 \times 3 \times 32$ array.
The output from each filter may be described as a \emph{feature map}.
The filters are applied with unit stride. In this way, each of the 32 feature maps has dimension 
$32 \times 32$. 
Max pooling is then applied to each feature map using stride length two.
This reduces the dimension of the feature maps to 
$16 \times 16$. A ReLU activation is then used. 
\item[Block 2] applies 
convolution followed 
 by activation and then 
 a pooling layer.
This
reduces the dimension to 
 $ 8 \times 8 \times 32$.
In more detail, we use $32$ filters. Each is  $5 \times 5$ across the dimensions of the feature maps, 
and also scans across all $32$ feature maps. 
So the weights could be regarded as a $ 5 \times 5 \times 32 \times 32$ tensor. 
The stride length is one, so the resulting 
$32$ feature maps 
are still of dimension $16 \times 16$.
After ReLU activation,  
an average pooling layer of stride two is then applied, which reduces each of the 
$32$ feature maps to dimension $8 \times 8$.
 \item[Block 3] applies a convolution layer followed by the activation function, and then performs a pooling operation, in a way that reduces dimension  to  $ 4 \times 4 \times 64$.
In more detail, $64$ filters are applied. Each filter is $5 \times 5$
across the dimensions of the feature maps, 
and also scans across all $32$ feature maps.  
So the weights could be regarded as a $ 5 \times 5 \times 32 \times 64$ tensor. 
The stride has length one, resulting in feature maps of dimension $8 \times 8$.
After ReLU activation, 
an average pooling layer of stride two is applied, which reduces each of the 
$64$ feature maps to dimension $4 \times 4$.
 \item[Block 4] does not use pooling, just convolution followed by activation, leading to dimension 
 $ 1 \times 1 \times 64$.
In more detail, $64$ filters are used. Each filter is 
 $ 4 \times 4$ across the $64$ feature maps,  so
the weights could be regarded as a $ 4 \times 4 \times 64 \times 64$ tensor, and 
each filter produces a single number.
\item[Block 5] does not involve convolution. It uses a general  
(fully connected) weight matrix of the type discussed in 
sections~\ref{sec:annex}
to 
\ref{sec:mat} to give output of dimension 
 $ 1 \times 1 \times 10$. 
This corresponds to a weight matrix of dimension $10 \times 64$. 
\item[A final softmax] operation transforms each of the ten ouput components to the range $[0,1]$. 
\end{description}
Figure~\ref{Fig:netblocks} gives a visual overview of the network architecture.

Our output is a vector of ten real numbers.
The cost function in the optimization problem 
takes the softmax log loss form 
(\ref{eq:smlcost}) with $K = 10$.
We specify stochastic gradient \emph{with momentum}, which  
uses a \lq\lq moving average\rq\rq\ of current and past gradient directions.
We use mini-batches of size 100
(so $m = 100$ in (\ref{eq:sgdefbatch})) and set a fixed number of 
$45$ epochs. We predefine the learning rate for each epoch: 
$\eta = 0.05$,   
$\eta = 0.005$ and
$\eta = 0.0005$ for the first $30$ epochs, next $10$ epochs and final $5$ epochs, respectively.
Running on 
a Tesla C2075 GPU in single precision,
the 45 epochs can be completed in just under 4 hours.

As an additional test, we also train the network with dropout.
Here, on each stochastic gradient step, 
any neuron has its output re-set to zero with independent  
probability
\begin{itemize}
 \item $0.15$ in block 1,
 \item $0.15$ in block 2,
 \item $0.15$ in block 3,
 \item $0.35$ in block 4,
 \item $0$ in block 5 (no dropout).
\end{itemize}
We emphasize that in this case all neurons become active when the trained network is applied to the test data.

\begin{figure}
\begin{center}
\includegraphics[scale=0.6]{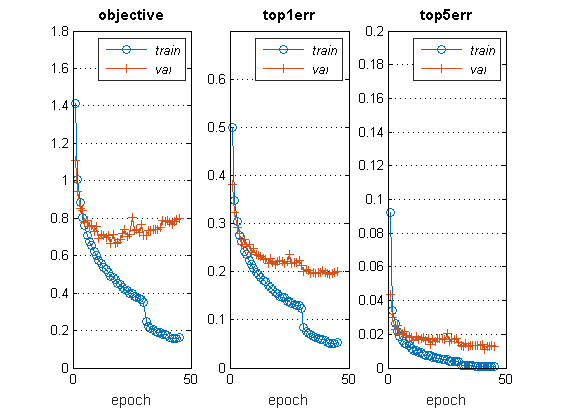}
\caption{
Errors for the trained network.
Horizontal axis runs over the 45 epochs of the stochastic gradient method (that is, 45 passes through the training data).
Left: Circles show cost function on the training data; crosses show cost function on the
validation data.
Middle:
Circles show the percentage of instances where the most likely classification from the network does not match the correct category, over the training data images; crosses show the same measure 
computed over the validation data.
Right:   
Circles show the percentage of instances where the five most 
likely classifications from the network do not include the 
correct category, over the training data images; 
crosses show the same measure computed over the validation data.
}
\label{Fig:neterr}
\end{center}
\end{figure}

In Figure~\ref{Fig:neterr} we illustrate the training process in the case of no dropout.
For the plot on the left, circles are used to show how the objective function 
(\ref{eq:smlcost}) 
decreases after each of the 45 epochs. 
We also use crosses to indicate the objective function value
on the validation data. 
(More precisely, these error measures are averaged over the individual batches that form the epoch---note that weights are updated after each batch.)
Given that our overall aim is to assign images to one of the ten classes, the middle plot 
in Figure~\ref{Fig:neterr} 
looks at the percentage of errors that take place 
when we classify with the highest probability choice. 
Similarly, the plot on the right shows the percentage of 
cases where the correct category is not among 
the top five.
We see from Figure~\ref{Fig:neterr} that the validation error starts to plateau at a stage where the 
stochastic gradient method continues to make significant reductions on the training error.
This gives an indication that we are overfitting---learning fine details about the training data that 
will not help the network to generalize to unseen data.
  
\begin{figure}
\begin{center}
\includegraphics[scale=0.6]{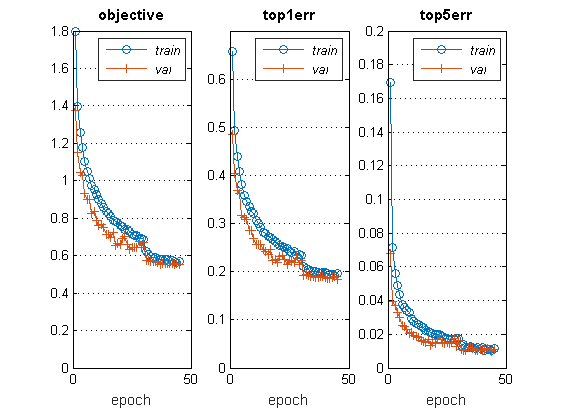}
\caption{
As for Figure~\ref{Fig:neterr} in the case where dropout was used.
}
\label{Fig:neterrdrop}
\end{center}
\end{figure}

Figure~\ref{Fig:neterrdrop}  shows the analogous results in the case where dropout is used. 
We see that the training errors are significantly 
larger than those in Figure~\ref{Fig:neterr} 
and the validation errors are of a similar magnitude. 
However, two 
key features in the dropout case are that (a) the validation error is below the training error,
and (b) the validation 
error continues to decrease in sync with the training error, both of which indicate that 
the optimization procedure is extracting useful information over all epochs. 

\begin{figure}
\begin{center}
\includegraphics[scale=0.6]{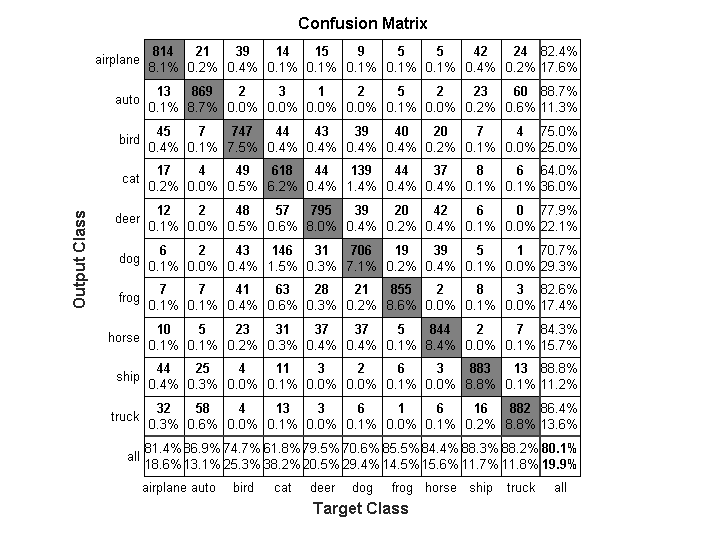}
\caption{
Confusion matrix for the the trained network from Figure~\ref{Fig:neterr}.  
}
\label{Fig:netconf}
\end{center}
\end{figure}

Figure~\ref{Fig:netconf} gives a summary of 
the performance of the trained network with no dropout (after 45 epochs)
in the form of a \emph{confusion matrix}.
Here, 
the integer value in the general $i,j$ entry shows the number of occasions where 
the network predicted category $i$ for an image from category
$j$.
  Hence, off-diagonal elements indicate mis-classifications.
For example, the 
(1,1) element equal to 814 in Figure~\ref{Fig:netconf} records the number of airplane images that were correctly classified as airplanes, and  
the 
(1,2) element equal to 21 records the number of automobile images that 
were incorrectly classified as airplanes.
Below each integer is the corresponding percentage, rounded to one decimal place, given that 
the test data has 1,000 images from each category.
The extra row, labeled \lq\lq all\rq\rq, summarizes the entries in each column. For example, the 
value 81.4\%\ in the first column of the final row arises because 814 of the 1,000 airplane images were correctly classified. Beneath this, the value 18.6\%\ arises because 186 of these airplane
images were incorrectly classified. 
The final column of the matrix, also labeled \lq\lq all\rq\rq, summarizes each row. For example, the
value 82.4\%\  in the final column of the first row arises because 988 images were classified by the 
network as airplanes, with 814 of these classifications being correct.
Beneath this, the value 17.6\%\ arises because the 
remaining 174 out of these 988
airplane classifications were incorrect.
Finally, the entries in the lower right corner summarize over all categories.
We see that 80.1\%\ of all images were correctly classified (and hence 
19.9\% were incorrectly classified). 

\begin{figure}
\begin{center}
\includegraphics[scale=0.6]{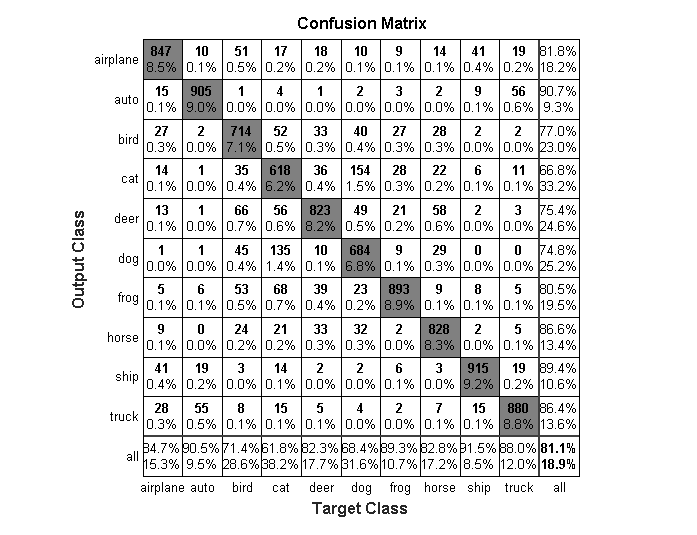}
\caption{
Confusion matrix for the the trained network from Figure~\ref{Fig:neterrdrop}, which used dropout.  
}
\label{Fig:netconfdrop}
\end{center}
\end{figure}

Figure~\ref{Fig:netconfdrop} gives the corresponding results in the case where dropout was used.
We see that the use of dropout has generally improved performance, and 
in particular has increased the overall success rate 
from 80.1\%\ to 81.1\%. Dropout gives larger values along the diagonal elements of the confusion matrix in nine out of the ten categories.

\begin{figure}
\begin{center}
\includegraphics[scale=0.7]{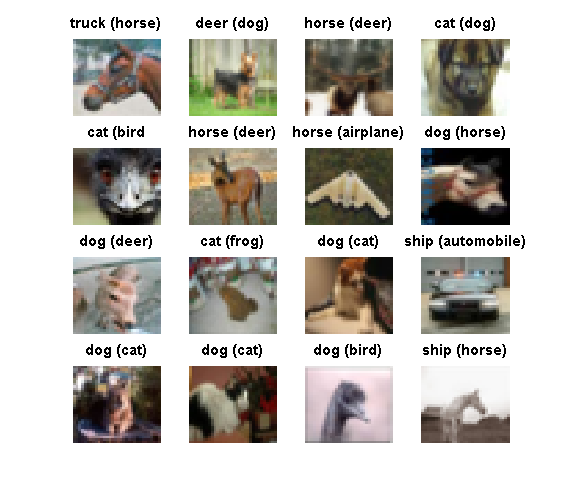}
\caption{
Sixteen of the images that were misclassified by the trained network from Figure~\ref{Fig:neterr}.
Predicted category is indicated, with correct category shown in parentheses. 
Note that images are low-resolution, having 
$32 \times 32$ pixels.
}
\label{Fig:netmiss}
\end{center}
\end{figure}

To give a feel for the difficulty of this task, Figure~\ref{Fig:netmiss} shows
16 images randomly sampled from those 
that were misclassified by the non-dropout network.

\section{Of Things Not Treated} \label{sec:disc}
This short introductory article is aimed at those who are new 
to deep learning.
In the interests of brevity and 
accessibility we have 
ruthlessly omitted 
many topics.
For those wishing to learn more, a good starting point 
is the 
free online book 
\cite{N15}, which provides a hands-on tutorial style 
description of deep learning techniques. 
The survey 
\cite{CBG15}
gives an intuitive and accessible 
overview of many of the key ideas behind deep learning, and highlights 
recent success stories. 
A more detailed overview of the
prize-winning performances
of deep learning tools can be found in
\cite{Schm15}, which also traces the development of
ideas across more than 800 references.
The review
\cite{WR17}
discusses the pre-history of deep learning and explains how
key ideas evolved.
For a comprehensive treatment of the state-of-the-art, 
we recommend the book 
\cite{GBC16} which, in particular, 
does an excellent job of 
introducing fundamental ideas from 
computer science/discrete mathematics, 
applied/computational mathematics and 
probability/statistics/inference before pulling them all  
together in the deep learning setting.
The recent review article 
\cite{BCN17}
focuses on optimization tasks arising in machine learning.
It summarizes the current theory underlying the 
stochastic gradient method, along with 
many 
alternative techniques. Those authors also emphasize that
optimization tools must be 
interpreted and judged carefully when operating 
within this inherently statistical framework.
Leaving aside the training issue, 
a mathematical framework for understanding the cascade 
of linear and nonlinear transformations used by deep 
networks is given in \cite{M16}.

To give a feel for some of the key issues that can be followed up, we finish with  
a list of questions that may have occured to 
interested readers, along with brief answers and 
further citations.

\begin{description}
\item[Why use artificial neural networks?]
Looking at 
Figure~\ref{Fig:pic_bdy}, it is clear that there are 
many ways to come up with a mapping that divides the x-y axis 
into two regions; a 
shaded region containing the circles and an 
unshaded region containing the crosses.
Artificial neural networks provide one useful approach.
In real applications, success corresponds to a small
\emph{generalization error}; the mapping should 
perform well when presented with new data.
In order to make  
rigorous, general, statements about performance,
we need to make some assumptions about 
the type of data. For example, we could 
analyze the situation where the data  
consists of samples drawn independently 
from a  certain probability distribution. If an 
algorithm is trained on such data, how will it perform 
when presented with \emph{new data from the same distribution}?
The authors in \cite{HRS16} prove that 
artificial neural networks trained with 
the stochastic gradient method can behave well in this sense.
Of course, in practice we 
cannot rely on the existence of such a distribution. 
Indeed, experiments in 
\cite{ZBHRV17} 
indicate that the worst case can be as bad as possible.
These authors 
tested state-of-the-art convolutional 
networks for image classification. 
In terms of the heuristic performance indicators used to monitor
the progress of the training phase, they found that  
the stochastic gradient method 
appears to work just as effectively when the 
\emph{images are randomly re-labelled}. This implies that the 
network is happy to learn noise---if the labels for the unseen data are similarly 
randomized then the classifications from the trained network are no better than random choice.  
Other authors have established negative results by showing that 
small and seemingly unimportant perturbations to an image can change its predicted class,
including cases where one pixel is altered \cite{SVK17}. 
Related work in \cite{BMRAG17} showed proof-of-principle for an  
\emph{adversarial patch}, which alters the classification when added to a wide range of images; for example, such a patch could be 
printed as a small sticker and used in the physical world.  
Hence, 
although artificial neural networks 
have outperformed rival methods in many application fields, 
the reasons behind this success 
are not fully understood. The survey 
\cite{vidalGBS17}
describes a range of mathematical approaches that
are beginning to provide useful insights, whilst the discussion piece
\cite{M18} includes a list of ten concerns.
\item[Which nonlinearity?]
The sigmoid function
(\ref{eq:sig}), illustrated in 
Figure~\ref{Fig:sigpic},
and the rectified linear unit (\ref{eq:relu}) are 
popular choices for the activation function.
Alternatives include
the \emph{step function},
\[
       \left\{ \begin{array}{c}
               0, \quad \mathrm{for~} x \le 0, \\
               1, \quad \mathrm{for~} x > 0. 
               \end{array}
                 \right.
\]
Each of these can 
undergo \emph{saturation}: produce very small 
derivatives that thereby 
reduce the size of the gradient updates.  
Indeed, 
the step function 
and 
 rectified linear unit have completely flat portions.
For this reason, 
a
 \emph{leaky rectified linear unit}, such as, 
\[
 f(x) = \left\{ \begin{array}{c}
               0.01x, \quad \mathrm{for~} x \le 0, \\
               x, \quad~~~~~ \mathrm{for~} x > 0, 
               \end{array}
                 \right.
\]
is sometimes preferred,  
in order to force a nonzero derivative for negative inputs.
The back propagation algorithm described in section~\ref{sec:bp} 
carries through to general activation functions.
\item[How do we decide on the structure of our net?]
Often, 
there is a natural choice for the size of the output layer. 
For example, to classify images of 
individual handwritten digits, it would make sense to have
an output layer consisting of ten 
neurons, corresponding to $0,1,2,\ldots,9$, 
as used in Chapter~1 of \cite{N15}.
In some cases, a physical application imposes natural constraints on one or more
of the hidden layers \cite{HMSPE18}. 
However, in general, choosing the overall number of layers, 
the number of neurons within each 
layer, and any constraints involving inter-neuron connections, 
is not an exact science. 
Rules of thumb have been suggested, but there is no widely 
accepted technique. 
In the context of image processing, it may be possible to attribute roles to different layers; 
for example, detecting edges, motifs and larger structures as information flows forward 
\cite{CBG15}, and our understanding of biological neurons provides further insights
\cite{GBC16}. 
But specific roles cannot be completely hardwired into the network design---the weights and biases, and hence the tasks performed by each layer, emerge from the training procedure.  
We note that the use of back propagation to 
compute gradients is not restricted 
to the types of connectivity, activation functions and cost functions 
discussed here. Indeed, the method fits into a very general 
framework of techniques known as 
\emph{automatic differentiation}  
or
\emph{algorithmic differentiation}  
\cite{GW08}.

\item[How big do deep learning networks get?]
The 
AlexNet
architecture \cite{AlexNet}
achieved groundbreaking image classification results in 2012.
This network used 650,000 neurons, with 
five convolutional layers followed by two fully connected layers and a final
softmax.
The programme \emph{AlphaGo}, developed by the Google 
DeepMind team to play the board game Go, rose to fame 
by beating 
the human European champion by five games to nil 
in October 2015 
\cite{silver16}.
AlphaGo 
makes use of 
two artificial neural networks with 
 13 layers and 15 layers,
some
convolutional and others fully connected, 
involving 
millions of weights. 
\item[Didn't my numerical analysis teacher tell me never to use
steepest descent?]
It is known that the steepest descent method can perform poorly on
examples where other methods, notably those using information about
the second derivative of the objective function,
are much more efficient. 
Hence, optimization textbooks typically downplay 
steepest descent \cite{Fl87,NW06}. 
However, it is important to note that
training an artificial neural network
is a very specific optimization task:
 \begin{itemize}
   \item the problem dimension and the expense of computing the
objective function and its derivatives, can be extremely high,
 \item the optimization task is set within a framework that is
inherently statistical in nature,
 \item a great deal of research effort has been devoted to the
development of practical improvements to the basic stochastic gradient
method in the deep learning context.
\end{itemize}
Currently, 
a theoretical underpinning for the success of the stochastic gradient method 
in training networks is far from complete \cite{BCN17}. 
A promising line of research is to connect 
the stochastic gradient method with discretizations of 
stochastic differential equations,  
\cite{SK17}, 
generalizing the idea that many deterministic optimization methods
can be viewed as timestepping methods for gradient ODEs, 
\cite{H99}.
  We also note that the introduction of more traditional tools
from the field of optimization may lead to improved 
training algorithms.

\item[Is it possible to regularize?]
As we discussed in section~\ref{sec:image}, 
overfitting  occurs when a trained network 
performs  accurately on the given data, but 
cannot \emph{generalize} well to new data.
\emph{Regularization} is a broad term that describes attempts to
avoid overfitting by rewarding smoothness.
One approach is to alter the cost function in order to encourage small weights.
For example, (\ref{eq:costfun}) could be extended to
\begin{equation}
\mathrm{Cost} 
     = 
      \frac{1}{N}
             \sum_{i=1}^{N}
                \hf \| 
              y( x^{\{ i \}} ) - a^{[L]}( x^{\{ i \}} )
                               \|_2^2
    +
 \frac{\lambda}{N}
   \sum_{l = 2}^{L}   \| \, W^{[l]} \, \|_2^2.
 \label{eq:costfunreg}
\end{equation}
Here $\lambda > 0$ is the regularization parameter. 
One motivation for 
(\ref{eq:costfunreg}) is that large weights may lead to neurons that are sensitive 
to their inputs, and hence less reliable when new data is presented. This argument does not 
apply to 
the biases, which typically are not included in such a regularization term.
It is straightforward to check that using 
 (\ref{eq:costfunreg})
instead of 
 (\ref{eq:costfun})
makes a very minor and inexpensive change to the back propagation algorithm.

\item[What about ethics and accountability?]
The use of \lq\lq algorithms\rq\rq\ to aid decision-making
is not a recent phenomenon. However, 
the increasing influence of black-box technologies
is naturally causing concerns in many quarters.
The recent articles
\cite{Dav17,Gr16}
raise several relevant issues and illustrate them with concrete examples.
They also highlight the 
particular challenges arising from massively-parameterized 
artificial neural networks.
Professional and governmental institutions
are, of course, alert to these matters.
In 2017, 
the Association for Computing Machinery's   
US Public Policy Council released seven 
\emph{Principles for Algorithmic Transparency and Accountability}\footnote{
\texttt{https://www.acm.org/} }. Among their recommendations  
are that 
\begin{itemize}
 \item \lq\lq  Systems and institutions that use algorithmic decision-making are encouraged to produce explanations regarding both the procedures followed
by the algorithm and the specific decisions that are made\rq\rq, 
and 
\item 
\lq\lq A description of the way in which the training data 
was collected should be maintained by the builders of 
the algorithms, 
accompanied by an exploration of the potential biases induced by the human or algorithmic data-gathering process.\rq\rq
\end{itemize}
Article 15 of the 
European Union's General Data Protection Regulation 
2016/679\footnote{\texttt{https://www.privacy-regulation.eu/en/15.htm}},
which takes effect in May 2018,
concerns \lq\lq Right of access by the data subject,\rq\rq\ and 
includes the requirement
that 
\lq\lq The data subject shall have the right to obtain from 
the controller confirmation as to whether or not personal data 
concerning him or her are being processed, and, where that is 
the case, access to the personal data and the following 
information:.\rq\rq Item (h) on the subsequent list covers 
\begin{itemize}
\item 
\lq\lq 
the existence of automated decision-making, including profiling, 
referred to in Article 22(1) and (4) and, at least in those cases, 
meaningful information about the logic involved, 
as well as the significance and the envisaged consequences 
of such processing for the data subject.\rq\rq
\end{itemize}

\item[What are some current research topics?]
Deep learning is 
a 
fast-moving, high-bandwith field,  
where many new advances are 
driven by the needs of specific 
application areas and the features of 
new high performance computing architectures.
Here, we briefly mention three hot-topic areas that 
have not yet been discussed.

Training a network can be an extremely expensive task.
When a trained network is seen to make a mistake on new data, it is therefore
tempting to fix this with a local perturbation to the weights and/or network structure, rather than 
re-training from scratch.
Approaches for 
this type of \emph{on the fly} tuning can be developed and justified using the 
theory of 
measure concentration in high dimensional spaces \cite{GT17}.

\emph{Adversarial networks}, 
\cite{GAN14},
are based on the 
concept that  an
artificial neural network 
may be viewed as a \emph{generative model}: a way to create
realistic data. Such a model may be useful, for example, 
as a means to produce realistic sentences, or very high resolution images. In the adversarial setting, the generative model is 
pitted against a 
\emph{discriminative model}.
The role of the 
discriminative model is to 
distinguish between real training data and data produced 
by the generative model.
By iteratively improving the performance of these models, the
quality of both the generation and discrimination can be 
increased dramatically.

The idea behind \emph{autoencoders} \cite{RHW86} is, perhaps surprisingly, to 
produce an overall network whose output matches its input.
More precisely, one network, 
known as the \emph{encoder}, 
corresponds to a map $F$ that takes an input vector, $x \in \RR^{s}$, and 
produces a lower dimensional output vector $F(x) \in \RR^{t}$.
So $ t \ll s$. 
Then a second network, 
known as the \emph{decoder}, 
corresponds to a map 
$G$ that takes us back to the same dimension as $x$; that is,
$G(F(x)) \in \RR^{s}$.
We could then aim to minimize the sum of the squared error
$\| x - G(F(x)) \|_2^2$ over a set of training data.
Note that this technique does not require the use of labelled 
data---in the case of images 
we are attempting to reproduce each picture without knowing 
what it depicts.
Intuitively, a good encoder is a tool for 
dimension reduction. It extracts the key features.
Similarly, a good decoder can reconstruct the data 
from those key features.

\item[Where can I find code and data?] 

There are many publicly available codes that provide access to 
deep learning algorithms.
In addition to 
{\sc{MatConvNet}}
\cite{Vedaldi15}, we mention 
Caffe \cite{jia2014caffe}, 
Keras \cite{Keras}, 
TensorFlow
\cite{45381},
Theano 
\cite{2016arXiv160502688full}  and 
Torch
\cite{Torch7}.
These packages differ in their underlying platforms and in the 
extent of expert knowledge 
required. 
Your favorite scientific computing 
environment 
may also offer 
a range of proprietary and user-contributed deep learning 
toolboxes.
However, it is currently the case that 
making serious use of 
modern deep learning technology 
requires a strong background in numerical computing.
Among the standard benchmark data sets 
are the 
CIFAR-10 collection
\cite{Krizhevsky09learningmultiple} that we used in 
section~\ref{sec:image}, and its big sibling CIFAR-100, 
ImageNet \cite{ImageNet}, 
and the 
handwritten digit database 
MNIST
\cite{lecun-98}.
\end{description}

\section*{Acknowledgements}
We are grateful
to the 
{\sc{MatConvNet}}
team 
for making their package available under 
a permissive BSD license. 
The MATLAB code in Listings~\ref{List.netbp}
and \ref{List.activate}  
can be found at 
\begin{verbatim}
http://personal.strath.ac.uk/d.j.higham/algfiles.html
\end{verbatim}
as well 
as an exteneded version that produces
Figures~\ref{Fig:pic_bp_cost} and 
\ref{Fig:pic_bdy_bp}, and
a 
MATLAB code that uses 
\texttt{lsqnonlin}
to 
produce
Figure~\ref{Fig:pic_bdy}.

\bibliographystyle{siam}
\bibliography{refsDL}

\end{document}